\documentclass[11pt,reqno]{amsart}
\pdfoutput=1

\usepackage[utf8]{inputenc}
\usepackage[T1]{fontenc}
\IfFileExists{lmodern.sty}{\usepackage{lmodern}}{}
\usepackage[a4paper,margin=1.15in]{geometry}

\usepackage{amsmath,amssymb,amsfonts,amsthm}
\usepackage{mathtools}
\usepackage{booktabs}
\usepackage{array}
\usepackage{enumitem}
\usepackage[final,expansion=false]{microtype}
\usepackage{xurl}
\usepackage{xcolor}
\usepackage[colorlinks=true,
linkcolor={blue!70!black},
citecolor={blue!70!black},
urlcolor={blue!80!black}]{hyperref}

\hypersetup{%
	pdftitle={Splitting types of Richelot isogenies over finite fields and a newform of weight four},
	pdfauthor={Hung Dang}}

\newtheorem{theorem}{Theorem}[section]
\newtheorem{proposition}[theorem]{Proposition}
\newtheorem{lemma}[theorem]{Lemma}
\newtheorem{corollary}[theorem]{Corollary}
\theoremstyle{definition}

\theoremstyle{remark}
\newtheorem{remark}[theorem]{Remark}

\DeclareMathOperator{\Res}{Res}
\DeclareMathOperator{\Disc}{Disc}
\DeclareMathOperator{\Jac}{Jac}
\DeclareMathOperator{\Stab}{Stab}
\providecommand{\PGL}{\operatorname{PGL}}
\newcommand{\F}{\mathbb{F}}
\newcommand{\Fq}{\mathbb{F}_q}
\newcommand{\Z}{\mathbb{Z}}
\newcommand{\Ql}{\mathbb{Q}_\ell}
\newcommand{\PP}{\mathbb{P}}
\newcommand{\dcoef}{\mathfrak{d}}
\newcommand{\Adm}{\mathcal{A}}
\newcommand{\cS}{\mathcal{S}}
\newcommand{\cT}{\mathcal{T}}
\newcommand{\cK}{\mathcal{K}}

\begin{document}

\title[Splitting types of Richelot isogenies and a newform of weight four]{Splitting types of Richelot isogenies over finite fields and a newform of weight four}

\author{Hung Dang}
\address{Department of Mathematics, University of Phuong Dong, Hanoi, Vietnam}
\email{hung.dt@phuongdong.edu.vn}
\thanks{ORCID 0009-0006-3272-0573.}


\begin{abstract}
	A factorization \(f=uvw\) of a square-free sextic over a finite field \(\Fq\)
	of odd characteristic into three monic quadratics defines a genus-2 curve
	with a rational \((2,2)\)-isogeny kernel, and the Richelot construction
	attaches to it three dual quadratics \(U,V,W\). We determine, for every odd
	prime power \(q\), the number of such factorizations with prescribed
	splitting behaviour over \(\Fq\) of all six quadratics. Each of the \(64\)
	counts is \(q(q-1)(q+1-n_f)/64\) times \(C(q)\pm t_q\), where \(n_f\) is the
	number of roots of \(f\) in \(\Fq\), \(C\) is an explicit monic cubic, and
	\(t_q\) is the trace of Frobenius attached to the newform
	\(\eta(2z)^{4}\eta(4z)^{4}\) of weight \(4\) and level \(8\). The proof reduces
	the counts by elementary means to a twisted second moment of the Legendre
	family of elliptic curves, which is a value of a hypergeometric function over
	\(\Fq\) evaluated by Ahlgren and Ono for primes. As consequences, the number of
	factorizations whose three pairwise resultants lie in one square class is a
	polynomial in \(q\), with proportion \(1/4+3/(4q^{2})+O(q^{-4})\). When the
	quotient is a product of two elliptic curves conjugate over \(\F_{q^{2}}\), we
	prove that the splitting of the factors \(u,v,w\) over \(\Fq\) decides which
	points of order two of the elliptic factor are divisible by two.
\end{abstract}

\maketitle

\section{Introduction}\label{sec:intro}

Let \(q\) be an odd prime power. A square-free sextic \(f\in\Fq[x]\) written as a
product \(f=uvw\) of three monic quadratics defines a genus-2 curve
\(C\colon y^{2}=f(x)\) together with a maximal isotropic subgroup of
\(\Jac(C)[2]\) defined over \(\Fq\), the kernel of a \((2,2)\)-isogeny. The
classical construction of Richelot describes the quotient. With the brackets
\[
U=v'w-vw',\qquad V=w'u-wu',\qquad W=u'v-uv',
\]
and \(\dcoef\) the determinant of the coefficient matrix of \((u,v,w)\), the
quotient is the Jacobian of \(\dcoef\,y^{2}=UVW\) when \(\dcoef\neq0\), and a
product of two elliptic curves when \(\dcoef=0\)
\cite{CasselsFlynn1996,Smith2005}. The roots of \(U\) are the fixed points of the
involution of \(\PP^{1}\) that exchanges the two roots of \(v\) and the two roots
of \(w\), and \(\Disc(U)=4\Res(v,w)\) when \(U\) has degree two; so \(U\) splits
over \(\Fq\) exactly when \(\Res(v,w)\) is a square (a bracket of degree one
counting as split), and cyclically.

This paper determines how the splitting behaviour of the three quadratics
\(u,v,w\) and that of the three dual quadratics \(U,V,W\) are distributed and
correlated over \(\Fq\). The answer is exact for every \(q\), and it involves a
modular form.

\subsection*{The 64 splitting types}
Let \(S\) be the set of the \(q(q-1)\) separable monic quadratics over \(\Fq\). An
ordered triple \((u,v,w)\in S^{3}\) is \emph{admissible} if \(u,v,w\) are pairwise
coprime, and \(\Adm(\Fq)\) denotes the set of admissible triples. Let \(\psi_q\) be
the quadratic character of \(\Fq\), with \(\psi_q(0)=0\). For an admissible triple
put
\begin{gather*}
e_s=\psi_q(\Disc s)\quad(s\in\{u,v,w\}),\\
r_u=\psi_q(\Res(v,w)),\qquad r_v=\psi_q(\Res(w,u)),\qquad r_w=\psi_q(\Res(u,v)).
\end{gather*}
Thus \(e_s=1\) when \(s\) splits over \(\Fq\) and \(e_s=-1\) when it is
irreducible, and \(r_u=1\) exactly when the two fixed points of the involution
attached to \((v,w)\) are rational, that is, when the bracket \(U\) splits (a
bracket of degree one counting as split, with a root at infinity). For
\(e,r\in\{\pm1\}^{3}\) let \(N(e;r)\) be the number of admissible triples with
\((e_u,e_v,e_w)=e\) and \((r_u,r_v,r_w)=r\).

Let \(g(z)=\eta(2z)^{4}\eta(4z)^{4}=\sum_{n\ge1}a(n)e^{2\pi inz}\) be the newform
of weight \(4\) for \(\Gamma_0(8)\), so that \(a(3)=-4\), \(a(5)=-2\), \(a(7)=24\).
For an odd prime \(p\) let \(\alpha_p,\beta_p\) be the roots of
\(X^{2}-a(p)X+p^{3}\), and put
\begin{equation}\label{eq:tq}
	t_{p^{n}}=\alpha_p^{n}+\beta_p^{n}\qquad(n\ge1).
\end{equation}
Thus \(t_q\) is the trace of the \(q\)-power Frobenius on the two-dimensional
\(\ell\)-adic representation attached to \(g\). It equals \(a(q)\) when \(q\) is
prime but not at proper prime powers; for instance \(t_9=a(3)^{2}-2\cdot3^{3}=-38\),
whereas \(a(9)=-11\).

\begin{theorem}\label{thm:A}
	Let \(q\) be an odd prime power, let \(e,r\in\{\pm1\}^{3}\), let \(j\) be the
	number of indices \(i\) with \(e_i=-1\), and put
	\(\sigma=(-1)^{j+1}r_1r_2r_3\). Then
	\begin{equation}\label{eq:A}
		N(e;r)=\frac{q(q-1)(q-5+2j)}{64}\,\bigl(C_{e;r}(q)+\sigma\,t_q\bigr),
	\end{equation}
	where \(C_{e;r}\) is the monic cubic of Table~\ref{tab:classes}.
\end{theorem}

\begin{table}[ht]
	\centering
	\caption{The cubics of Theorem~\ref{thm:A}, grouped into the \(20\) orbits of
		the simultaneous permutations of \((u,v,w)\) and \((U,V,W)\). Here \(n_{ab}\) is
		the number of \(s\in\{u,v,w\}\) with \((e_s,r_s)=(a,b)\); for instance
		\(n_{+-}\) counts the \(s\) that split while the opposite bracket is
		irreducible. The column \(\#\) is the size of the orbit, and
		\(j=n_{-+}+n_{--}\).}
	\label{tab:classes}
	\small
	\begin{tabular}{@{}ccccccll@{}}
		\toprule
		\(n_{++}\) & \(n_{+-}\) & \(n_{-+}\) & \(n_{--}\) & \(\#\) & \(j\) & \(C_{e;r}(q)\) & \(\sigma\)\\
		\midrule
		3 & 0 & 0 & 0 & 1 & 0 & \(q^{3}-12q^{2}+56q-91\) & \(-1\) \\
		2 & 1 & 0 & 0 & 3 & 0 & \(q^{3}-10q^{2}+30q-23\) & \(+1\) \\
		1 & 2 & 0 & 0 & 3 & 0 & \(q^{3}-8q^{2}+16q-7\) & \(-1\) \\
		0 & 3 & 0 & 0 & 1 & 0 & \(q^{3}-6q^{2}+14q-11\) & \(+1\) \\
		\midrule
		2 & 0 & 1 & 0 & 3 & 1 & \(q^{3}-6q^{2}+14q-11\) & \(+1\) \\
		1 & 1 & 1 & 0 & 6 & 1 & \((q+1)(q^{2}-5q+5)\) & \(-1\) \\
		2 & 0 & 0 & 1 & 3 & 1 & \(q^{3}-4q^{2}+8q-3\) & \(-1\) \\
		0 & 2 & 1 & 0 & 3 & 1 & \((q+1)(q^{2}-3q+1)\) & \(+1\) \\
		1 & 1 & 0 & 1 & 6 & 1 & \((q+1)(q^{2}-3q+1)\) & \(+1\) \\
		0 & 2 & 0 & 1 & 3 & 1 & \((q+1)(q^{2}-q+1)\) & \(-1\) \\
		\midrule
		1 & 0 & 2 & 0 & 3 & 2 & \((q+1)(q^{2}-5q+5)\) & \(-1\) \\
		0 & 1 & 2 & 0 & 3 & 2 & \((q+1)(q^{2}-3q+1)\) & \(+1\) \\
		1 & 0 & 1 & 1 & 6 & 2 & \((q+1)(q^{2}-3q+1)\) & \(+1\) \\
		0 & 1 & 1 & 1 & 6 & 2 & \((q+1)(q^{2}-q+1)\) & \(-1\) \\
		1 & 0 & 0 & 2 & 3 & 2 & \((q+1)(q^{2}-q-7)\) & \(-1\) \\
		0 & 1 & 0 & 2 & 3 & 2 & \((q+1)(q^{2}+q-3)\) & \(+1\) \\
		\midrule
		0 & 0 & 3 & 0 & 1 & 3 & \(q^{3}-6q^{2}-2q+37\) & \(+1\) \\
		0 & 0 & 2 & 1 & 3 & 3 & \((q+1)(q^{2}-5q+5)\) & \(-1\) \\
		0 & 0 & 1 & 2 & 3 & 3 & \((q+1)(q^{2}-3q+1)\) & \(+1\) \\
		0 & 0 & 0 & 3 & 1 & 3 & \(q^{3}-8q+9\) & \(-1\) \\
		\bottomrule
	\end{tabular}
\end{table}

The factor \(q-5+2j\) equals \(q+1-n_f\), where \(n_f=6-2j\) is the number of
roots of \(f=uvw\) in \(\Fq\), and when the class is not empty
\(q(q-1)(q+1-n_f)\) is the number of elements of \(\PGL_2(\Fq)\) that move no root
to infinity. Accordingly \eqref{eq:A} does not
depend on the affine chart: \(\tfrac1{64}\bigl(C_{e;r}(q)+\sigma t_q\bigr)\) is the
number of \(\PGL_2(\Fq)\)-orbits, weighted by the inverse orders of the
stabilisers, of ordered triples of pairwise disjoint reduced divisors of degree
two on \(\PP^{1}\) over \(\Fq\) of type \((e;r)\)
(Proposition~\ref{prop:intrinsic}(b)). The term \(t_q\) enters every one of the
\(64\) classes with coefficient \(\pm1\), and it cancels when one sums over the
types of \(u,v,w\) or over those of \(U,V,W\).

\begin{corollary}\label{cor:B}
	Let \(q\) be an odd prime power.
	\begin{enumerate}[label=\textup{(\alph*)},leftmargin=*]
		\item The number of admissible triples for which \(u,v,w,U,V,W\) all split
		over \(\Fq\) is
		\[
		\frac{q(q-1)(q-5)\bigl(q^{3}-12q^{2}+56q-91-t_q\bigr)}{64}.
		\]
		\item The number of admissible triples whose three resultants are squares,
		that is, for which \(U,V,W\) all split, is
		\(\tfrac18\,q(q-1)(q-3)(q^{3}-5q^{2}+12q-24)\).
		\item The number of admissible triples whose three resultants lie in one
		square class of \(\Fq^{\times}\), that is, for which \(U,V,W\) are all split
		or all irreducible, is
		\[
		A_q^{=}=\frac{q(q-1)\bigl(q^{4}-5q^{3}+15q^{2}-35q+40\bigr)}{4},
		\qquad
		\frac{A_q^{=}}{\#\Adm(\Fq)}=\frac14+\frac{3}{4q^{2}}+O(q^{-4}).
		\]
	\end{enumerate}
\end{corollary}

For example \(t_7=24\), and (a) gives \(42\) triples over \(\F_7\). Part (c) says
that the three square classes agree as often as three independent fair signs
would, up to \(O(q^{-2})\). The two outcomes of the construction deviate from
\(1/4\) in opposite directions: when \(\dcoef=0\) the three resultants always lie
in one class (Proposition~\ref{prop:pencil}), and among the triples with
\(\dcoef\neq0\) the proportion is \(1/4-3/(4q)+O(q^{-3})\)
(Remark~\ref{rem:quarter}).

\subsection*{The decomposable case}
When \(\dcoef=0\) the three quadratics lie in one pencil, and the common square
class of the resultants decides whether the two elliptic factors of the quotient
are defined over \(\Fq\) or are conjugate over \(\F_{q^{2}}\). In the second case
\(\Jac(C)\) has the Weil polynomial of the Weil restriction of an elliptic curve
\(E_1\) over \(\F_{q^{2}}\), to which it is therefore isogenous over \(\Fq\) by
Tate's theorem \cite[\S3, Thm.~1(c)]{Tate1966}, and
\(\#\Jac(C)(\Fq)=\#E_1(\F_{q^{2}})\); see Proposition~\ref{prop:factors}. These
facts are known: for the covers see \cite[Prop.~1]{CastryckDecruSmith2020}, and
for the field of definition of the elliptic factors and the Weil restriction see
\cite[Prop.~3.4]{BruinDoerksen2011} and \cite{BernsteinLange2014}. In this case
the splitting of the brackets carries no information, and we show that the
splitting of \(u,v,w\) themselves controls the \(2\)-primary part of
\(E_1(\F_{q^{2}})\).

\begin{theorem}\label{thm:C}
	Let \((u,v,w)\) be an admissible triple over \(\Fq\) with \(\dcoef=0\) whose
	resultants are nonsquares, let \(E_1\) be the elliptic curve \eqref{eq:E1} over
	\(K=\F_{q^{2}}\), and let \(k\) be the number of the factors \(u,v,w\) that split
	over \(\Fq\), so that \(k=3-j\) with \(j\) as in Theorem~\ref{thm:A}. Then \(E_1\) has full \(2\)-torsion over \(K\), with points
	\(P_u,P_v,P_w\) of order two indexed by the factors, and \(P_u\in2E_1(K)\) if and
	only if \(v\) and \(w\) both split over \(\Fq\), and cyclically. Consequently,
	with \(N=\#\Jac(C)(\Fq)=\#E_1(K)\):
	\begin{enumerate}[label=\textup{(\alph*)},leftmargin=*]
		\item if \(k=3\), then \(E_1[4]\subseteq E_1(K)\), so \(16\mid N\);
		\item if \(k=2\), then the \(2\)-primary part of \(E_1(K)\) is
		\(\Z/2\times\Z/2^{m}\) with \(m\ge2\), so \(8\mid N\) and
		\(E_1[4]\not\subseteq E_1(K)\);
		\item if \(k\le1\), then the \(2\)-primary part of \(E_1(K)\) is \((\Z/2)^{2}\),
		so \(N\equiv4\pmod 8\).
	\end{enumerate}
\end{theorem}

Jacobians of this kind are used for genus-2 arithmetic with elliptic point
counting, and Bernstein and Lange report that in their construction, where all
six Weierstrass points are rational, computer experiments suggest that the
subgroup of order \(16\) of the elliptic curve over \(\F_{p^{2}}\) is
\((\Z/4)^{2}\) \cite[\S5]{BernsteinLange2014}; part (a) gives this when their
elliptic curve is \(E_1\) or its conjugate (Remark~\ref{rem:BL}).

\subsection*{Method and relation to earlier work}
The proof of Theorem~\ref{thm:A} is by Walsh--Fourier inversion on
\(\{\pm1\}^{6}\). The \(64\) sums of the monomials in the six signs over
\(\Adm(\Fq)\) are computed from a lemma on signed sums over Frobenius-stable
pairs of points and from point counts on curves of genus at most one
(\S\ref{subsec:lemmas}). Sixty of them are polynomials in \(q\)
(Table~\ref{tab:walsh}). The remaining four, those in which the product
\(r_ur_vr_w\) is weighted by at least two of the \(e_s\), reduce to the
\emph{twisted second moment of the Legendre family}
\begin{equation}\label{eq:Tintro}
	\cT(q)=\sum_{\lambda\in\Fq\setminus\{0,1\}}\psi_q(\lambda)\,a_\lambda(q)^{2},
\end{equation}
where \(a_\lambda(q)\) is the trace of Frobenius of \(y^{2}=x(x-1)(x-\lambda)\)
(Proposition~\ref{prop:reduction}). The second moment arises because the signed
sum over the third quadratic of a triple is the sum of the traces of Frobenius
on \(H^{0}\), \(H^{1}\) and the symmetric square of \(H^{1}\) of the genus-one
curve \(y^{2}=u(x)v(x)\)
(Remark~\ref{rem:geometry}). One has \(\cT(q)+1=q^{3}\,{}_4F_3(1)_q\), a value of
Greene's hypergeometric function over \(\Fq\), and Ahlgren and Ono proved that
\(p^{3}\,{}_4F_3(1)_p=-a(p)-p\) for every odd prime \(p\)
\cite[Thm.~6]{AhlgrenOno2000}; an equivalent identity follows from the trace
formulas of \cite[Prop.~2.1(2)]{FrechetteOnoPapanikolas2004}, which those
authors attribute to another paper of Ahlgren and Ono. For prime powers we
deduce \(\cT(q)=-t_q-q-1\) from the prime case by a standard \(\ell\)-adic
argument (Theorem~\ref{thm:T}); this step is implicit in the cohomological
description of rigid Calabi--Yau threefolds (cf.\
\cite[Thm.~4.1 and Table~1]{SaitoYui2001}), and moments of this kind are treated
for all prime powers by the trace formula in \cite{KaplanPetrow2017}.

The modular form and the threefolds that carry it are known. Van Geemen and
Nygaard proved that \(g\) gives the \(L\)-function of \(H^{3}\) of a Siegel
modular threefold \cite[Thm.~2.4]{vanGeemenNygaard1995}. The double cover of
\(\PP^{3}\) branched along
\(x_0x_1x_2x_3(x_0-x_1)(x_1-x_2)(x_2-x_3)(x_3-x_0)=0\), of which the affine
threefold \(\mathcal X\) of Proposition~\ref{prop:hyp}, with
\(q^{3}+1+\cT(q)\) points, is a chart, is the double octic of type
\(T^{1}_{70}\) in the list of Cynk and Meyer of rigid
Calabi--Yau threefolds attached to \(g\) \cite[\S2.2]{CynkMeyer2005}. Cynk,
Freitag and Salvati Manni study the Siegel threefold of the group
\(\Gamma_2[2]\cap\Gamma_{2,0}[4]\) and state that a rigid Calabi--Yau model of it
is modular with \(g\) \cite[\S4]{CynkFreitagSalvatiManni2011}. For the class of
Corollary~\ref{cor:B}(a), in which all six quadratics split, the cofactor of
the count is the number of \(\Fq\)-points of an open subscheme
\(\mathcal Y^{\circ}\) of the model of \cite{CynkFreitagSalvatiManni2011},
defined over \(\Z[1/2]\) by their equations and the nonvanishing of ten
explicit forms. Indeed,
the threefold \(\mathcal R\), whose \(\Fq\)-points are the classes of all-split
Richelot data with ordered roots and ordered fixed points, is isomorphic over
\(\Z[1/2]\) to \(\mathcal Y^{\circ}\) (Proposition~\ref{prop:cfs} and
Corollary~\ref{cor:cfs}). That \(\mathcal Y^{\circ}\) is the open part of the
Siegel modular threefold on which the ten fourth powers of the even theta
constants do not vanish rests on the description of the ring of modular forms
in \cite{CynkFreitagSalvatiManni2011}, whose proof is sketched there, and on the
choice of model (Remark~\ref{rem:cfs-model}). For this class the reduction to \(\cT(q)\) is short
(Remark~\ref{rem:geometry}). What we
add is the reduction from the Richelot data
for all \(64\) classes and the resulting table, in which \(t_q\) occurs with a
sign depending on the class; we are not aware of a statement of this table, or
of the polynomial count of Corollary~\ref{cor:B}(c), in the works cited or in the
literature on point counts of moduli of genus-2 curves with level structure that
we consulted. Since \(\sigma=-e_1e_2e_3r_1r_2r_3\), in the projective
normalisation \(t_q\) enters only through the character \(e_1e_2e_3r_1r_2r_3\),
the product of all six signs. We have no proof of the corresponding statement
about the cohomology of the threefold, and we did not find that statement in the
works cited; we formulate it as a question in Remark~\ref{rem:cfs-signs}.

By a classical formula for translations by points of order two on the Kummer
surface, the triples of Corollary~\ref{cor:B}(c) are those whose kernel acts
diagonally on rational Kummer coordinates. We use that formula without
re-deriving it and therefore state this as a remark (Remark~\ref{rem:kummer}),
which nothing else in the paper uses.

Richelot isogenies are the degree-two steps of explicit isogeny computations in
genus 2, in particular in isogeny-based cryptography, where the decomposable case
is the event that a codomain splits. The results here concern arithmetic
statistics over a general odd \(\Fq\), not such computations. The counts of the admissible triples by outcome of the
construction (decomposable, or quotient curve of degree five or six), the
corresponding masses, and a classification of the outcome from the minors of the
coefficient matrix are the subject of the companion paper
\cite{DangNguyen2026v2}; the present paper does not depend on it.

\subsection*{Organization}
Section~\ref{sec:prelim} recalls the Richelot data and the decomposable case.
Section~\ref{sec:walsh} computes the \(64\) Walsh sums, evaluates the moment
\(\cT(q)\) and proves Theorem~\ref{thm:A} and Corollary~\ref{cor:B};
Subsection~\ref{subsec:cfs} identifies the all-split stratum with the open
subscheme \(\mathcal Y^{\circ}\) of the model of Cynk, Freitag and Salvati Manni.
Section~\ref{sec:dec} proves Theorem~\ref{thm:C}. Section~\ref{sec:comp} lists the
exact computations that accompany the proofs. It names the two places where
proofs depend on them, and the statements made in remarks that were checked only
by machine; nothing else in the paper uses those statements.

\section{Richelot data}\label{sec:prelim}

Throughout, \(q\) is an odd prime power, \(\psi_q\) is the quadratic character of
\(\Fq\) with \(\psi_q(0)=0\), and \(\psi_{q^{2}}\) is that of \(\F_{q^{2}}\). The set
\(S\) of separable monic quadratics over \(\Fq\) is the disjoint union of the set
\(S_{\mathrm{sp}}\) of those that split over \(\Fq\) and the set
\(S_{\mathrm{ir}}\) of the irreducible ones; each has \(q(q-1)/2\) elements. Roots
are taken in \(\F_{q^{2}}\).

\subsection{Minors, brackets and the involution of a pair}

For monic quadratics \(a=x^{2}+a_1x+a_0\) and \(b=x^{2}+b_1x+b_0\) over a field of
odd characteristic put
\begin{equation}\label{eq:minors}
	M_2=b_1-a_1,\qquad M_1=b_0-a_0,\qquad M_0=a_1b_0-a_0b_1,
\end{equation}
and \([a,b]=a'b-ab'\).

\begin{lemma}\label{lem:minors}
	With this notation:
	\begin{enumerate}[label=\textup{(\roman*)},leftmargin=*,itemsep=0pt]
		\item \([a,b]=M_2x^{2}+2M_1x+M_0\);
		\item \(M_0-a_1M_1+a_0M_2=0\) and \(M_0-b_1M_1+b_0M_2=0\);
		\item \(a(x)b(z)-a(z)b(x)=(x-z)\bigl(M_2xz+M_1(x+z)+M_0\bigr)\);
		\item \(\Res(a,b)=M_1^{2}-M_2M_0\), and hence \(\Disc([a,b])=4\Res(a,b)\) when
		\(M_2\neq0\).
	\end{enumerate}
\end{lemma}

\begin{proof}
	(i), (ii) and (iii) are verified by expanding. (iv) Since \(a\) is monic,
	\(\Res(a,b)=\Res(a,b-a)=\prod_{a(\varrho)=0}(M_2\varrho+M_1)
	=M_1^{2}-a_1M_1M_2+a_0M_2^{2}\), and \(a_1M_1-a_0M_2=M_0\) by (ii).
\end{proof}

Let \(a,b\) be coprime. By (iii), two points \(x\neq z\) of \(\PP^{1}\) have the
same image under \(x\mapsto[a(x):b(x)]\) exactly when
\(M_2xz+M_1(x+z)+M_0=0\), read in homogeneous form when one of them is
\(\infty\). So
\begin{equation}\label{eq:iota}
	\iota_{a,b}(x)=\frac{-M_1x-M_0}{M_2x+M_1}
\end{equation}
is the involution of \(\PP^{1}\) that exchanges the two roots of \(a\) and the two
roots of \(b\); its matrix has determinant \(-\Res(a,b)\neq0\). Its fixed points
are the roots of \([a,b]\), together with \(\infty\) when \(M_2=0\). They are
distinct, and by (iv) they are rational over \(\Fq\) if and only if \(\Res(a,b)\)
is a square in \(\Fq\): when \(M_2=0\) the fixed points are \(\infty\) and
\(-M_0/(2M_1)\), and \(\Res(a,b)=M_1^{2}\). We say that the bracket \([a,b]\)
\emph{splits} over \(\Fq\) in this case, a bracket of degree one counting as
split.

\subsection{Admissible triples}

An ordered triple \((u,v,w)\in S^{3}\) is \emph{admissible} if \(u,v,w\) are
pairwise coprime; equivalently, \(f=uvw\) is square-free of degree six. The set
of admissible triples is \(\Adm(\Fq)\). For such a triple let
\[
U=[v,w],\qquad V=[w,u],\qquad W=[u,v],\qquad
\dcoef=\det\begin{pmatrix}u_0&u_1&1\\v_0&v_1&1\\w_0&w_1&1\end{pmatrix},
\]
where \(s=x^{2}+s_1x+s_0\). The classes of the divisors of zeros of \(u,v,w\),
minus the divisor at infinity of \(C\colon y^{2}=f\), form with \(0\) a maximal
isotropic subgroup \(G\) of \(\Jac(C)[2]\) defined over \(\Fq\). If \(\dcoef\neq0\),
then \(UVW\) is square-free of degree \(5\) or \(6\) and \(\Jac(C)/G\) is the
Jacobian of \(\dcoef\,y^{2}=UVW\); if \(\dcoef=0\), the quotient is a product of two
elliptic curves \cite{CasselsFlynn1996}, \cite[Ch.~8]{Smith2005}. We call
an admissible triple with \(\dcoef=0\) \emph{decomposable}. The signs \(e_s\) and
\(r_s\) of the introduction are defined for every admissible triple, and
\(r_u=1\) exactly when \(U\) splits over \(\Fq\), by Lemma~\ref{lem:minors}(iv).

\subsection{The decomposable case}

\begin{proposition}\label{prop:pencil}
	Let \((u,v,w)\) be an admissible triple with \(\dcoef=0\).
	\begin{enumerate}[label=\textup{(\roman*)},leftmargin=*,itemsep=0pt]
		\item There are \(a,b\in\Fq^{\times}\) with \(a+b=1\) and \(w=au+bv\), and
		\(U=-aW\), \(V=-bW\). The three involutions \(\iota_{v,w}\), \(\iota_{w,u}\),
		\(\iota_{u,v}\) are equal; we write \(\iota\) for this involution.
		\item The resultants \(\Res(v,w)\), \(\Res(w,u)\), \(\Res(u,v)\) lie in one
		square class of \(\Fq^{\times}\), which we call the \emph{class} of the
		triple. If the class is a nonsquare, then \(M_2(u,v)\neq0\).
		\item No fixed point of \(\iota\) is a root of \(uvw\).
	\end{enumerate}
\end{proposition}

\begin{proof}
	(i) The rows of the matrix are dependent and any two of them are independent,
	so \(w=au+bv\); comparing leading coefficients gives \(a+b=1\), and \(ab\neq0\)
	because \(w\neq u,v\). Then \(U=v'(au+bv)-v(au'+bv')=-aW\), and \(V=-bW\) in the
	same way. The involution of a pair is determined by its fixed points, which
	are the zeros of the bracket on \(\PP^{1}\). (ii) By (i) and
	Lemma~\ref{lem:minors}, \(\Res(v,w)=a^{2}\Res(u,v)\) and
	\(\Res(w,u)=b^{2}\Res(u,v)\); if \(M_2(u,v)=0\), then \(\Res(u,v)=M_1^{2}\).
	(iii) Let \(\xi\) be a fixed point with \(u(\xi)=0\). Then
	\(0=W(\xi)=u'(\xi)v(\xi)\), which is impossible because \(u\) is separable and
	coprime to \(v\). The same argument applies to \(v\) with \(W\), and to \(w\)
	with \(U\).
\end{proof}

Let \((u,v,w)\) be decomposable with \(M_2=M_2(u,v)\neq0\), let \(\alpha,\beta\) be
the roots of \(W\), that is, the fixed points of \(\iota\), and let
\(K=\Fq(\alpha)\). Thus \(K=\Fq\) if the class is a square and \(K=\F_{q^{2}}\)
otherwise, and
\begin{equation}\label{eq:alphabeta}
	\alpha+\beta=-\frac{2M_1}{M_2},\qquad \alpha\beta=\frac{M_0}{M_2},\qquad
	(\beta-\alpha)^{2}=\frac{4\Res(u,v)}{M_2^{2}} .
\end{equation}

\begin{proposition}\label{prop:factors}
	With this notation put \(T=(x-\alpha)/(x-\beta)\), so that
	\(x=(\alpha-\beta T)/(1-T)\), and \(\tilde y=(1-T)^{3}y\).
	\begin{enumerate}[label=\textup{(\roman*)},leftmargin=*,itemsep=0pt]
		\item For \(s\in\{u,v,w\}\) one has
		\((1-T)^{2}\,s\bigl(x(T)\bigr)=s(\beta)\,T^{2}+s(\alpha)\), and
		\(s(\alpha)s(\beta)\neq0\). Hence \(C_K\) is the curve
		\(\tilde y^{2}=F(T^{2})\), where
		\begin{equation}\label{eq:E1}
			F(X)=\prod_{s\in\{u,v,w\}}\bigl(s(\beta)X+s(\alpha)\bigr),
			\qquad E_1\colon Y^{2}=F(X),\qquad E_2\colon Y^{2}=X^{3}F(1/X),
		\end{equation}
		and \(E_1,E_2\) are elliptic curves over \(K\) with the maps of degree two
		\(\pi_1\colon(T,\tilde y)\mapsto(T^{2},\tilde y)\) and
		\(\pi_2\colon(T,\tilde y)\mapsto(T^{-2},\tilde yT^{-3})\) from \(C_K\).
		\item \(\Jac(C_K)\) is isogenous to \(E_1\times E_2\) by
		\((\pi_{1,*},\pi_{2,*})\), with kernel \(G\).
		\item If the class is a nonsquare, then \(E_2\) is the conjugate of \(E_1\)
		under \(\mathrm{Gal}(\F_{q^{2}}/\Fq)\), and
		\(P_{\Jac(C)/\Fq}(Z)=P_{E_1/\F_{q^{2}}}(Z^{2})\), where \(P(Z)=\det(1-Z\,\Phi)\)
		is the characteristic polynomial of the Frobenius endomorphism \(\Phi\) on the
		\(\ell\)-adic Tate module. In particular \(\#\Jac(C)(\Fq)=\#E_1(\F_{q^{2}})\).
	\end{enumerate}
\end{proposition}

\begin{proof}
	(i) The coefficient of \(T\) in
	\((\alpha-\beta T)^{2}+s_1(\alpha-\beta T)(1-T)+s_0(1-T)^{2}\) is
	\(-\bigl(2\alpha\beta+s_1(\alpha+\beta)+2s_0\bigr)=-2(M_0-s_1M_1+s_0M_2)/M_2\).
	This vanishes for \(s=u\) and \(s=v\) by Lemma~\ref{lem:minors}(ii), and for
	\(s=w=au+bv\) because the expression is affine in \((s_1,s_0)\) and \(a+b=1\).
	The coefficients of \(T^{2}\) and \(T^{0}\) are \(s(\beta)\) and \(s(\alpha)\),
	which are nonzero by Proposition~\ref{prop:pencil}(iii). Since \(f\) is
	square-free, so are \(F\) and \(X^{3}F(1/X)\), and the formulas for \(\pi_1\) and
	\(\pi_2\) follow from the equations. (ii) The maps \(\pi_1,\pi_2\) are the
	quotients by the involution \((T,\tilde y)\mapsto(-T,\tilde y)\) and
	by its product with the hyperelliptic involution, so
	\((\pi_{1,*},\pi_{2,*})\) is an isogeny of degree four: by the Lemma of
	\cite[\S2]{Kuhn1988} the pull-back \(E_1\times E_2\to\Jac(C)\) is an isogeny
	with kernel of order four, and its composition with \((\pi_{1,*},\pi_{2,*})\) is
	multiplication by two.
	The two roots \(\pm\varrho\) of a factor \(s(\beta)T^{2}+s(\alpha)\) have the
	same image under \(\pi_1\) and under \(\pi_2\), which gives the kernel; compare
	\cite[Prop.~1 and \S3.2]{CastryckDecruSmith2020}. (iii) The nontrivial
	automorphism \(\tau\) of \(K/\Fq\) exchanges \(\alpha\) and \(\beta\), so it
	sends \(F\) to \(X^{3}F(1/X)\), \(T\) to \(1/T\) and \(\tilde y\) to
	\(-\tilde y\,T^{-3}\). Hence \(E_1^{\tau}=E_2\) and \(\pi_1^{\tau}=[-1]\circ\pi_2\).
	The \(q\)-power Frobenius morphisms therefore satisfy
	\(\mathrm{Fr}_{E_1}\circ\pi_1=[-1]\circ\pi_2\circ\mathrm{Fr}_{C}\), so the
	\(q\)-power Frobenius \(\Phi\) of \(\Jac(C)\) exchanges the two summands
	\(\pi_i^{*}V_\ell E_i\) of \(V_\ell\Jac(C)\). In an adapted basis it has the block
	form \(\left(\begin{smallmatrix}0&B\\A&0\end{smallmatrix}\right)\), so
	\(\det(1-Z\Phi)=\det(1-Z^{2}BA)\), and \(BA\) is the action of \(\Phi^{2}\) on
	\(\pi_1^{*}V_\ell E_1\), that is, of the \(q^{2}\)-power Frobenius of \(E_1\).
	Evaluating at \(Z=1\) gives the last assertion.
\end{proof}

The statements of Proposition~\ref{prop:factors} are known: for the covers see
\cite[Prop.~1]{CastryckDecruSmith2020}, and for the field of definition of the
elliptic factors and Jacobians isogenous to a Weil restriction see
\cite[Prop.~3.4]{BruinDoerksen2011} and \cite{BernsteinLange2014}. We include the short proof because
the explicit form \eqref{eq:E1} of \(E_1\) is used in \S\ref{sec:dec}. If
\(M_2(u,v)=0\), then \(u_1=v_1=w_1\), the class is a square, and after the
translation \(x\mapsto x-u_1/2\) the sextic is \(F_0(x^{2})\) with a cubic \(F_0\); the
two factors are then \(Y^{2}=F_0(X)\) and \(Y^{2}=XF_0(X)\) over \(\Fq\).

The decomposable triples are counted as follows. In the plane of the
coefficients \((s_1,s_0)\), a triple is decomposable exactly when its three
points are collinear.

\begin{lemma}\label{lem:Dq}
	The number of decomposable admissible triples is
	\(D_q=q(q-1)(q-2)(q^{2}-3q+3)\), and \(\tfrac12q^{2}(q-1)^{2}(q-2)\) of them
	have nonsquare class. Two quadratics on a line \(s_0=\mu s_1+\nu\) with middle
	coefficients \(s_1\neq s_1'\) have resultant \((s_1'-s_1)^{2}(\mu^{2}+\nu)\).
\end{lemma}

\begin{proof}
	On the line \(s_0=\mu s_1+\nu\) the minors of two points are \(M_2=s_1'-s_1\),
	\(M_1=\mu M_2\) and \(M_0=-\nu M_2\), which gives the resultant by
	Lemma~\ref{lem:minors}(iv). So the points of such a line are pairwise coprime
	exactly when \(\mu^{2}+\nu\neq0\). The discriminant restricted to the line is
	\(s_1^{2}-4\mu s_1-4\nu\), with discriminant \(16(\mu^{2}+\nu)\); so the line
	carries \(q\) separable quadratics if \(\mu^{2}+\nu\) is a nonsquare and \(q-2\)
	if it is a nonzero square, and there are \(q(q-1)/2\) lines of each kind. Each
	of the \(q\) lines \(s_1=\mathrm{const}\) carries \(q-1\) separable quadratics,
	any two with resultant \(M_1^{2}\neq0\). Counting ordered triples of distinct
	points,
	\[
	D_q=q(q-1)(q-2)(q-3)+\tfrac12q(q-1)\bigl[q(q-1)(q-2)+(q-2)(q-3)(q-4)\bigr],
	\]
	which is the stated polynomial, and the middle term is the number of triples
	of nonsquare class.
\end{proof}

The count \(D_q\) also appears in the companion paper \cite{DangNguyen2026v2};
the proof above is self-contained, and the present paper does not depend on that
paper.

\section{Splitting types and the newform}\label{sec:walsh}

This section proves Theorem~\ref{thm:A} and Corollary~\ref{cor:B}. The signs
\(e_s\) and \(r_s\) of an admissible triple are those of the introduction; for
\(s\in S\) we also write \(e_s=\psi_q(\Disc s)\), so that \(s\) has \(1+e_s\) roots
in \(\Fq\).

\subsection{Walsh sums}

For \(\varepsilon,\rho\in\{0,1\}^{3}\) define the \emph{Walsh sum}
\[
\cS[\varepsilon;\rho]=\sum_{(u,v,w)\in\Adm(\Fq)}
e_u^{\varepsilon_1}e_v^{\varepsilon_2}e_w^{\varepsilon_3}\,
r_u^{\rho_1}r_v^{\rho_2}r_w^{\rho_3}.
\]
Thus \(\cS[000;000]=\#\Adm(\Fq)\), and by inversion on \(\{\pm1\}^{6}\), with
\(e^{\varepsilon}=e_1^{\varepsilon_1}e_2^{\varepsilon_2}e_3^{\varepsilon_3}\) and
\(r^{\rho}=r_1^{\rho_1}r_2^{\rho_2}r_3^{\rho_3}\),
\begin{equation}\label{eq:inversion}
	N(e;r)=\frac1{64}\sum_{\varepsilon,\rho\in\{0,1\}^{3}}
	e^{\varepsilon}r^{\rho}\,\cS[\varepsilon;\rho].
\end{equation}
The set \(\Adm(\Fq)\) is stable under permutations of \((u,v,w)\), so \(N\) and
\(\cS\) do not change when one permutation is applied to the three indices of
both arguments. We write
\(|\varepsilon|=\varepsilon_1+\varepsilon_2+\varepsilon_3\).

For \(\lambda\in\Fq\) put
\(a_\lambda(q)=-\sum_{x\in\Fq}\psi_q\bigl(x(x-1)(x-\lambda)\bigr)\); for
\(\lambda\neq0,1\) this is the trace of Frobenius of the Legendre curve
\(y^{2}=x(x-1)(x-\lambda)\). With the twisted second moment \(\cT(q)\) of
\eqref{eq:Tintro} put
\begin{equation}\label{eq:T}
	\cK(q)=\cT(q)+q+1 .
\end{equation}

\begin{theorem}\label{thm:walsh}
	Let \(q\) be an odd prime power and let \(t_q\) be as in \eqref{eq:tq}. Then
	\(\cS[111;111]=-q(q-1)(q-2)\,t_q\), and
	\(\cS[\varepsilon;111]=q(q-1)(t_q+2)\) for \(|\varepsilon|=2\). The other \(60\)
	Walsh sums are the polynomials in \(q\) of Table~\ref{tab:walsh}.
\end{theorem}

	\begin{table}[ht]
		\centering
		\caption{The Walsh sums \(\cS[\varepsilon;\rho]\) divided by \(q(q-1)\).
			For \(|\rho|\in\{0,3\}\) put \(k=|\varepsilon|\). For \(|\rho|=1\) let \(i\) be
			the index with \(\rho_i=1\), for \(|\rho|=2\) the index with \(\rho_i=0\), and
			put \(m=|\varepsilon|-\varepsilon_i\). By Theorem~\ref{thm:T},
			\(\cK(q)=-t_q\).}
		\label{tab:walsh}
		\footnotesize\setlength{\tabcolsep}{4pt}
		\begin{tabular}{@{}cclll@{}}
			\toprule
			\(|\rho|\) & \(\varepsilon_i\) & & & \\
			\midrule
			\(0\) & & \(k=0\): \((q-2)^{2}(q^{2}-q+4)\) & \(k=1\): \(-(q-2)(2q^{2}-7q+7)\) & \\
			& & \(k=2\): \(-(q-2)(q^{2}-5q+8)\) & \(k=3\): \((q-2)(5q-7)\) & \\
			\(1\) & \(0\) & \(m=0\): \(-(q^{3}-4q^{2}+8q-9)\) & \(m=1\): \(2q^{2}-7q+7\) & \(m=2\): \(-(3q-7)\)\\
			\(1\) & \(1\) & \(m=0\): \((2q-3)(q-3)\) & \(m=1\): \(q^{2}-5q+8\) & \(m=2\): \(-(3q-5)\)\\
			\(2\) & \(0\) & \(m=0\): \(q^{2}-5q+8\) & \(m=1\): \((q-2)(q-3)\) & \(m=2\): \(-(3q-7)\)\\
			\(2\) & \(1\) & \(m=0\): \(-(3q-5)\) & \(m=1\): \(-(q-5)\) & \(m=2\): \(-(q-3)\)\\
			\(3\) & & \(k=0\): \(-(q-5)\) & \(k=1\): \(-(q-3)\) & \\
			& & \(k=2\): \(2-\cK(q)\) & \(k=3\): \((q-2)\cK(q)\) & \\
			\bottomrule
		\end{tabular}
	\end{table}

Propositions~\ref{prop:reduction} and~\ref{prop:poly} express the Walsh sums
through \(q\) and \(\cK(q)\) by elementary means, and Theorem~\ref{thm:T} shows
that \(\cK(q)=-t_q\).

	\subsection{Two lemmas on character sums}\label{subsec:lemmas}

	For \(h\in\Fq[X]\) put
	\(\Sigma_1(h)=\sum_{a\in\Fq}\psi_q(h(a))\),
	\(\Sigma_2(h)=\sum_{a\in\F_{q^{2}}}\psi_{q^{2}}(h(a))\) and
	\(n_h=\#\{a\in\Fq: h(a)=0\}\), where \(\psi_{q^{2}}\) is the quadratic character of
	\(\F_{q^{2}}\). For \(w\in S\) with roots \(a,b\in\F_{q^{2}}\) we have
	\(\Res(w,h)=h(a)h(b)\in\Fq\); this is multiplicative in \(h\).

	\begin{lemma}[Signed sums over Frobenius-stable pairs]\label{lem:A}
		For every \(h\in\Fq[X]\),
		\begin{align}
			\sum_{w\in S}e_w\,\psi_q\bigl(\Res(w,h)\bigr)&=\tfrac12\bigl(\Sigma_1(h)^{2}-\Sigma_2(h)\bigr),\label{eq:Aplus}\\
			\sum_{w\in S}\psi_q\bigl(\Res(w,h)\bigr)&=\tfrac12\bigl(\Sigma_1(h)^{2}+\Sigma_2(h)\bigr)-(q-n_h).\label{eq:Azero}
		\end{align}
	\end{lemma}

	\begin{proof}
		Taking roots identifies \(S_{\mathrm{sp}}\) with the unordered pairs of distinct
		elements of \(\Fq\), and \(S_{\mathrm{ir}}\) with the pairs \(\{a,a^{q}\}\),
		\(a\in\F_{q^{2}}\setminus\Fq\). So the part of \eqref{eq:Azero} over
		\(S_{\mathrm{sp}}\) is
		\(\tfrac12\sum_{a\neq b\in\Fq}\psi_q(h(a))\psi_q(h(b))=\tfrac12\bigl(\Sigma_1(h)^{2}-(q-n_h)\bigr)\).
		For \(z\in\F_{q^{2}}\) we have \(\psi_q(z^{q+1})=\psi_{q^{2}}(z)\), because
		\((z^{q+1})^{(q-1)/2}=z^{(q^{2}-1)/2}\), and \(h(a^{q})=h(a)^{q}\). So for
		\(w\in S_{\mathrm{ir}}\) with root \(a\) we get
		\(\psi_q(\Res(w,h))=\psi_{q^{2}}(h(a))\), and the part over \(S_{\mathrm{ir}}\) is
		\(\tfrac12\sum_{a\in\F_{q^{2}}\setminus\Fq}\psi_{q^{2}}(h(a))=\tfrac12\bigl(\Sigma_2(h)-(q-n_h)\bigr)\),
		since every element of \(\Fq^{\times}\) is a square in \(\F_{q^{2}}\). The sum of
		the two parts is \eqref{eq:Azero}, and their difference is
		\eqref{eq:Aplus}.
	\end{proof}

	\begin{lemma}[Curves of genus at most one]\label{lem:B}
		\begin{enumerate}[label=\textup{(\roman*)},leftmargin=*,itemsep=0pt]
			\item If \(h\in\Fq[X]\) is monic and squarefree of degree \(4\) and
			\(\theta(h)=-1-\Sigma_1(h)\), then \(\Sigma_2(h)=2q-1-\theta(h)^{2}\).
			\item If \(h\in\Fq[X]\) is squarefree of degree \(3\) and
			\(\theta(h)=-\Sigma_1(h)\), then \(\Sigma_2(h)=2q-\theta(h)^{2}\).
			\item If \(h=l\,(X^{2}+bX+c)\) with \(l\in\Fq^{\times}\) and \(b^{2}\neq4c\), then
			\(\Sigma_1(h)=-\psi_q(l)\) and \(\Sigma_2(h)=-1\).
			\item If \(h=(X-r)^{2}h_0\) with \(r\in\Fq\) and \(h_0\in\Fq[X]\), then
			\(\Sigma_1(h)=\Sigma_1(h_0)-\psi_q(h_0(r))\) and
			\(\Sigma_2(h)=\Sigma_2(h_0)-\psi_{q^{2}}(h_0(r))\).
		\end{enumerate}
	\end{lemma}

	\begin{proof}
		(i) The smooth projective model \(C_h\) of \(y^{2}=h(x)\) has genus one and two
		points at infinity, both rational because \(h\) is monic. So
		\(\#C_h(\F_{q^{n}})=q^{n}+2+\sum_{x\in\F_{q^{n}}}\psi_{q^{n}}(h(x))\), and
		\(\#C_h(\F_{q^{n}})=q^{n}+1-(\alpha^{n}+\beta^{n})\) with \(\alpha\beta=q\)
		\cite[Thm.~V.2.3.1]{Silverman2009}. The case \(n=1\) gives
		\(\alpha+\beta=\theta(h)\), and \(n=2\) gives
		\(\Sigma_2(h)=-1-(\theta(h)^{2}-2q)\). (ii) is the same with one point at
		infinity. (iii) is \cite[Thm.~5.48]{LidlNiederreiter1997} over \(\Fq\) and over
		\(\F_{q^{2}}\), where \(\psi_{q^{2}}(l)=1\). (iv) holds because
		\(\psi_q(h(a))=\psi_q(h_0(a))\) for \(a\neq r\) and \(\psi_q(h(r))=0\), and
		likewise over \(\F_{q^{2}}\).
	\end{proof}

	\subsection{The sums with \texorpdfstring{\(\rho=111\)}{rho=111}}

	For \(u,v\in S\) write \(R_{uv}=\Res(u,v)\) and
	\(\theta_{uv}=-1-\sum_{x\in\Fq}\psi_q\bigl(u(x)v(x)\bigr)\). Thus
	\(\theta_{uv}=\theta_{vu}\), and when \(u,v\) are coprime \(\theta_{uv}=\theta(uv)\)
	is the trace of Frobenius of the genus-one curve \(y^{2}=u(x)v(x)\). For
	\(a,b\in\{0,1\}\) and \(k\ge0\) put
	\[
	J_k^{ab}=\sum_{u,v\in S}e_u^{a}e_v^{b}\,\psi_q(R_{uv})\,\theta_{uv}^{k},
	\]
	which is symmetric in \((a,b)\); the factor \(\psi_q(R_{uv})\) vanishes unless
	\(u\) and \(v\) are coprime.

	\begin{proposition}\label{prop:reduction}
		For every odd prime power \(q\), the sum \(\cS[\varepsilon;111]\) depends
		only on \(k=|\varepsilon|\), and for \(k=0,1,2,3\) it equals respectively
		\[
		-q(q-1)(q-5),\quad -q(q-1)(q-3),\quad q(q-1)\bigl(2-\cK(q)\bigr),\quad
		q(q-1)(q-2)\,\cK(q).
		\]
	\end{proposition}

	\begin{proof}
		\emph{Step 1.} We first sum over the third quadratic of a triple. Let
		\(u,v\in S\) be coprime and
		\(\theta=\theta_{uv}\). Then
		\begin{equation}\label{eq:step1}
			\sum_{w\in S}e_w\,\psi_q(R_{uw}R_{vw})=\theta^{2}+\theta+1-q,\qquad
			\sum_{w\in S}\psi_q(R_{uw}R_{vw})=\theta+2+e_u+e_v .
		\end{equation}
		Indeed \(R_{uw}R_{vw}=\Res(w,uv)\), the polynomial \(h=uv\) is monic and
		squarefree of degree four with \(n_h=2+e_u+e_v\), and \(\Sigma_1(h)=-1-\theta\),
		\(\Sigma_2(h)=2q-1-\theta^{2}\) by Lemma~\ref{lem:B}(i); insert this in
		Lemma~\ref{lem:A}.

		\emph{Step 2.} Next we evaluate the sums over pairs. For \(u\in S\),
		Lemma~\ref{lem:A} with \(h=u\) and Lemma~\ref{lem:B}(iii), which gives
		\(\Sigma_1(u)=\Sigma_2(u)=-1\), yield
		\begin{equation}\label{eq:step2}
			\sum_{v\in S}e_v\,\psi_q(R_{uv})=1,\qquad \sum_{v\in S}\psi_q(R_{uv})=e_u+1-q .
		\end{equation}
		Since \(\sum_{u\in S}e_u=0\) and \(\sum_{u\in S}e_u^{2}=q(q-1)\), it follows that
		\(J_0^{00}=-q(q-1)^{2}\), \(J_0^{10}=J_0^{01}=q(q-1)\) and \(J_0^{11}=0\).

		\emph{Step 3.} We evaluate sums over one quadratic weighted by its values
		at given points. Let \(u\in S\) and \(x\in\Fq\) with
		\(u(x)\neq0\), and put \(h_{u,x}=(x-X)\,u\), a squarefree cubic with \(2+e_u\)
		roots in \(\Fq\). For \(v\in S\) with roots \(c,d\) we have
		\(R_{uv}\,v(x)=u(c)u(d)(x-c)(x-d)=\Res(v,h_{u,x})\), so by
		Lemmas~\ref{lem:B}(ii) and~\ref{lem:A}
		\begin{equation}\label{eq:step3}
			\sum_{v\in S}e_v\,\psi_q\bigl(R_{uv}v(x)\bigr)=\theta(h_{u,x})^{2}-q,\qquad
			\sum_{v\in S}\psi_q\bigl(R_{uv}v(x)\bigr)=2+e_u .
		\end{equation}
		Similarly, let \(x,c,d\in\Fq\) with \(c\neq x\), \(d\neq x\), and
		\(\kappa=(x-X)(c-X)(d-X)\), so that \(u(x)u(c)u(d)=\Res(u,\kappa)\) for \(u\in S\). Then
		\begin{equation}\label{eq:step4}
			\sum_{u\in S}e_u\,\psi_q\bigl(u(x)u(c)u(d)\bigr)=\theta(\kappa)^{2}-q\ \ (c\neq d),
			\qquad =1\ \ (c=d).
		\end{equation}
		If \(c\neq d\), then \(\kappa\) is a squarefree cubic. If \(c=d\),
		Lemma~\ref{lem:B}(iv) with \(h_0=x-X\) gives \(\Sigma_1(\kappa)=-\psi_q(x-c)\)
		and \(\Sigma_2(\kappa)=-1\), and \eqref{eq:Aplus} gives \((1+1)/2\).

		\emph{Step 4.} We compute the first moments \(J_1^{ab}\), namely
		\begin{equation}\label{eq:step5}
			J_1^{00}=q(q-1)(q+1),\qquad J_1^{10}=J_1^{01}=0,\qquad J_1^{11}=-q(q-1)\,\cK(q).
		\end{equation}
		By definition
		\(\psi_q(R_{uv})\,\theta_{uv}=-\psi_q(R_{uv})-\sum_{x}\psi_q(u(x))\,\psi_q(R_{uv}v(x))\).
		The factor \(\psi_q(u(x))\) removes the \(x\) with \(u(x)=0\), so
		\eqref{eq:step3} applies to every remaining term, and
		\(\sum_x\psi_q(u(x))=-1\) by Lemma~\ref{lem:B}(iii). Summing over \(v\) and
		then over \(u\) gives \(J_1^{00}=-J_0^{00}+\sum_{u}(2+e_u)\),
		\(J_1^{10}=-J_0^{10}+\sum_{u}e_u(2+e_u)=0\), and \(J_1^{11}=-\Omega\) with
		\(\Omega=\sum_{u}e_u\sum_{x}\psi_q(u(x))\,\theta(h_{u,x})^{2}\), because
		\(J_0^{11}=0\) and \(\sum_ue_u\sum_x\psi_q(u(x))=0\). Expanding the square,
		\[
		\Omega=\sum_{x,c,d\in\Fq}\psi_q\bigl((x-c)(x-d)\bigr)\sum_{u\in S}e_u\,\psi_q\bigl(u(x)u(c)u(d)\bigr).
		\]
		The terms with \(c=x\) or \(d=x\) vanish. By \eqref{eq:step4} the terms
		with \(c=d\neq x\) contribute \(q(q-1)\), and those with \(x,c,d\) distinct
		contribute \(\sum\psi_q((x-c)(x-d))(\theta(\kappa)^{2}-q)\). For distinct \(x,c,d\)
		put \(\lambda=(d-x)/(c-x)\in\Fq\setminus\{0,1\}\). The substitution
		\(X=x+(c-x)X'\) gives \(\kappa=-(c-x)^{3}X'(X'-1)(X'-\lambda)\), hence
		\(\theta(\kappa)^{2}=a_\lambda(q)^{2}\), and \((x-c)(x-d)=(c-x)^{2}\lambda\). The map
		\((x,c,d)\mapsto(x,c-x,\lambda)\) is a bijection from the triples of distinct
		elements onto \(\Fq\times\Fq^{\times}\times(\Fq\setminus\{0,1\})\), and
		\(\sum_{\lambda\neq0,1}\psi_q(\lambda)=-1\), so
		\[
		\Omega=q(q-1)+q(q-1)\sum_{\lambda\neq0,1}\psi_q(\lambda)\bigl(a_\lambda(q)^{2}-q\bigr)
		=q(q-1)\,\cK(q).
		\]

		\emph{Step 5.} We treat the case \(|\varepsilon|\le2\). The summand
		\(e^{\varepsilon}\psi_q(R_{uv})\,\psi_q(R_{uw}R_{vw})\) of
		\(\cS[\varepsilon;111]\) vanishes unless the triple is admissible, so the sum may be taken over \(S^{3}\), and it depends
		only on \(|\varepsilon|\) because the weight \(\psi_q(R_{uv}R_{vw}R_{wu})\) is
		symmetric. Summing over \(w\) by the second formula of \eqref{eq:step1}
		and using \(e_u^{2}=1\),
		\begin{align*}
			\cS[000;111]&=J_1^{00}+2J_0^{00}+2J_0^{10}=-q(q-1)(q-5),\\
			\cS[100;111]&=J_1^{10}+2J_0^{10}+J_0^{00}+J_0^{11}=-q(q-1)(q-3),\\
			\cS[110;111]&=J_1^{11}+2J_0^{11}+2J_0^{10}=q(q-1)\bigl(2-\cK(q)\bigr).
		\end{align*}

		\emph{Step 6.} We reduce the case \(|\varepsilon|=3\) to a second
		moment. By the first formula of \eqref{eq:step1},
		\(\cS[111;111]=J_2^{11}+J_1^{11}+(1-q)J_0^{11}=J_2^{11}+J_1^{11}\). Put
		\(Z_{uv}=(1+\theta_{uv})^{2}\) and
		\(J_*=\sum_{u,v\in S}e_ue_v\,\psi_q(R_{uv})\,Z_{uv}\), so that
		\(J_2^{11}=J_*-2J_1^{11}-J_0^{11}\) and \(\cS[111;111]=J_*-J_1^{11}\). It
		remains to show that \(J_*=q(q-1)(q-3)\,\cK(q)\). Since
		\(Z_{uv}=\sum_{x,y\in\Fq}\psi_q(u(x)u(y))\,\psi_q(v(x)v(y))\),
		\begin{equation}\label{eq:Jexp}
			J_*=\sum_{u\in S}e_u\sum_{x,y\in\Fq}\psi_q(u(x)u(y))\sum_{v\in S}e_v\,\psi_q\bigl(R_{uv}v(x)v(y)\bigr).
		\end{equation}
		For \(x\neq y\) in \(\Fq\) let
		\(v_{xy}=(X-x)(X-y)\in S_{\mathrm{sp}}\). For \(u\in S\) and \(x,y\in\Fq\) with
		\(u(x)u(y)\neq0\) we have \(R_{uv}v(x)v(y)=\Res(v,h)\) with
		\(h=u\cdot(X-x)(X-y)\), and \eqref{eq:Aplus} gives
		\begin{equation}\label{eq:7ab}
			\sum_{v\in S}e_v\,\psi_q\bigl(R_{uv}v(x)v(y)\bigr)=
			\begin{cases}
				\theta'^{2}+\theta'+1-q,\quad\theta'=\theta_{u,v_{xy}}, & x\neq y,\\
				2+\psi_q(u(x)), & x=y.
			\end{cases}
		\end{equation}
		For \(x\neq y\) this is Step~1, since \(h=u\,v_{xy}\). For \(x=y\),
		Lemma~\ref{lem:B}(iii),(iv) give \(\Sigma_1(h)=-1-\psi_q(u(x))\) and
		\(\Sigma_2(h)=-2\). In \eqref{eq:Jexp} the terms with \(x=y\) therefore
		contribute
		\(\sum_ue_u\bigl(2(q-1-e_u)-1\bigr)=-2q(q-1)\). The ordered pairs \(x\neq y\)
		correspond two-to-one to the elements \(v_{xy}\) of \(S_{\mathrm{sp}}\), and
		\(\psi_q(u(x)u(y))=\psi_q(R_{u,v_{xy}})\). Since \(\theta^{2}+\theta+1=Z-\theta\),
		\begin{equation}\label{eq:Jstar}
			J_*=-2q(q-1)+2\sum_{v\in S_{\mathrm{sp}}}\Theta_v,\qquad
			\Theta_v:=\sum_{u\in S}e_u\,\psi_q(R_{uv})\bigl(Z_{uv}-\theta_{uv}-q\bigr).
		\end{equation}
		Fix \(v\in S_{\mathrm{sp}}\) with roots \(x,y\). Since \(R_{uv}=u(x)u(y)\),
		\[
		\sum_{u\in S}e_u\,\psi_q(R_{uv})\,Z_{uv}
		=\sum_{c,d\in\Fq}\psi_q(v(c)v(d))\sum_{u\in S}e_u\,\psi_q\bigl(R_{uv}u(c)u(d)\bigr),
		\]
		and only \(c,d\notin\{x,y\}\) contribute. The inner sum is \eqref{eq:7ab}
		with the roles of \(u\) and \(v\) exchanged. For \(c\neq d\) it is
		\(\theta^{2}+\theta+1-q\) with \(\theta=\theta_{v,v_{cd}}\), and
		\(\psi_q(v(c)v(d))=\psi_q(R_{v,v_{cd}})\); the ordered pairs \((c,d)\) correspond
		two-to-one to the \(v'\in S_{\mathrm{sp}}\) coprime to \(v\). For \(c=d\) it is
		\(2+\psi_q(v(c))\), and these terms add up to
		\(2(q-2)+\sum_c\psi_q(v(c))=2q-5\). Hence
		\begin{equation}\label{eq:7c}
			\sum_{u\in S}e_u\,\psi_q(R_{uv})\,Z_{uv}
			=2q-5+2\sum_{v'\in S_{\mathrm{sp}}}\psi_q(R_{vv'})\bigl(\theta_{vv'}^{2}+\theta_{vv'}+1-q\bigr).
		\end{equation}
		For signs \(\pm\) standing for \(S_{\mathrm{sp}}\) and \(S_{\mathrm{ir}}\) let
		\(J_k^{\pm\pm}\) be the sum of \(\psi_q(R_{uv})\theta_{uv}^{k}\) over \(u\) in the
		first and \(v\) in the second of the two sets indicated. Summing
		\eqref{eq:7c} over \(v\in S_{\mathrm{sp}}\) and using the first formula of
		\eqref{eq:step2},
		\begin{equation}\label{eq:7sum}
			\sum_{v\in S_{\mathrm{sp}}}\Theta_v=\tfrac12q(q-1)(q-5)
			+2\bigl(J_2^{++}+J_1^{++}+(1-q)J_0^{++}\bigr)-\bigl(J_1^{++}-J_1^{-+}\bigr).
		\end{equation}
		The three kinds of terms are evaluated as follows.

		(i) \(J_0^{++}=-q(q-1)(q-3)/4\), because for \(v\in S_{\mathrm{sp}}\) the sum of
		\(\psi_q(R_{vv'})\) over \(v'\in S_{\mathrm{sp}}\) is
		\(\tfrac12\sum_{c\neq d}\psi_q(v(c))\psi_q(v(d))=\tfrac12\bigl(1-(q-2)\bigr)\).

		(ii) By symmetry \(J_1^{+-}=J_1^{-+}\), and \(J_1^{00}=J_1^{++}+2J_1^{+-}+J_1^{--}\),
		\(J_1^{11}=J_1^{++}-2J_1^{+-}+J_1^{--}\), \(J_1^{10}=J_1^{++}-J_1^{--}=0\). Hence
		\(J_1^{++}=(J_1^{00}+J_1^{11})/4\) and \(J_1^{++}-J_1^{-+}=J_1^{11}/2\).

		(iii) \(J_2^{++}=q(q-1)(q-3)\,\cT(q)/4\). A coprime pair
		\((v_{xy},v_{cd})\) arises from exactly four quadruples \((x,y,c,d)\) of distinct
		elements of \(\Fq\), so \(J_2^{++}\) is \(\tfrac14\) of the sum of
		\(\psi_q\bigl((c-x)(c-y)(d-x)(d-y)\bigr)\,\theta(x,y,c,d)^{2}\) over these
		quadruples, where
		\(\theta(x,y,c,d)=-1-\sum_{a\in\Fq}\psi_q\bigl((a-x)(a-y)(a-c)(a-d)\bigr)\). The
		substitution \(a=x+(y-x)a'\) multiplies the quartic by \((y-x)^{4}\); with
		\(c'=(c-x)/(y-x)\) and \(d'=(d-x)/(y-x)\) it gives
		\(\theta(x,y,c,d)=\theta(0,1,c',d')\), and the character becomes
		\(\psi_q(c'(c'-1)d'(d'-1))\). As \((x,y,c,d)\mapsto(x,y-x,c',d')\) is a bijection
		onto \(\Fq\times\Fq^{\times}\times\{(c',d'):c',d'\notin\{0,1\},\ c'\neq d'\}\),
		\[
		J_2^{++}=\frac{q(q-1)}{4}\sum_{\substack{c,d\in\Fq\setminus\{0,1\}\\ c\neq d}}
		\psi_q\bigl(c(c-1)d(d-1)\bigr)\,\theta(0,1,c,d)^{2}.
		\]
		Fix such a pair \((c,d)\) and put \(\lambda=c(1-d)/(d(1-c))\). For \(a\neq0\)
		and \(s=1/a\) we have \(a(a-1)(a-c)(a-d)=s^{-4}(1-s)(1-cs)(1-ds)\), and the term
		\(s=0\) of \(\sum_{s\in\Fq}\psi_q((1-s)(1-cs)(1-ds))\) is \(1\); hence
		\(\theta(0,1,c,d)=-\sum_{s\in\Fq}\psi_q\bigl((1-s)(1-cs)(1-ds)\bigr)\). The
		substitution \(s=1+(1/c-1)s'\) turns \((1-s)(1-cs)(1-ds)\) into
		\(-cd\,(1/c-1)^{3}\,s'(s'-1)(s'-\lambda)\), so
		\(\theta(0,1,c,d)^{2}=a_\lambda(q)^{2}\). Here \(\lambda\neq0,1\) because
		\(c\neq0\), \(d\neq1\), \(c\neq d\), and \(c(c-1)d(d-1)=\lambda\,d^{2}(1-c)^{2}\),
		so the character is \(\psi_q(\lambda)\). Finally every
		\(\lambda\in\Fq\setminus\{0,1\}\) has exactly \(q-3\) preimages \((c,d)\): the
		equation is \(c=d\,(\lambda+c(1-\lambda))\); it has no solution for
		\(c=\lambda/(\lambda-1)\), and for each of the other \(q-3\) values of
		\(c\in\Fq\setminus\{0,1\}\) its unique solution \(d=c/(\lambda+c(1-\lambda))\)
		is different from \(0\), from \(1\) (otherwise \(\lambda(1-c)=0\)) and from
		\(c\) (otherwise \((1-\lambda)(c-1)=0\)).

		Inserting (i)--(iii) and \eqref{eq:step5} in \eqref{eq:7sum}, the
		two terms \(J_1^{11}/2\) cancel and
		\[
		\sum_{v\in S_{\mathrm{sp}}}\Theta_v
		=\tfrac12q(q-1)\bigl[(q-5)+(q-3)\cT(q)+(q+1)+(q-1)(q-3)\bigr],
		\]
		so that \(J_*=q(q-1)\bigl[(q-3)\cT(q)+q^{2}-2q-3\bigr]=q(q-1)(q-3)\,\cK(q)\)
		by \eqref{eq:Jstar}, and
		\(\cS[111;111]=J_*-J_1^{11}=q(q-1)(q-2)\,\cK(q)\).
	\end{proof}

\subsection{The polynomial components}

\begin{proposition}\label{prop:poly}
	For every odd prime power \(q\), the \(56\) Walsh sums \(\cS[\varepsilon;\rho]\)
	with \(\rho\neq111\) are the polynomials of Table~\ref{tab:walsh}.
\end{proposition}

The proof occupies the rest of this subsection. For \(u,v\in S\) the indicator of
\(\gcd(u,v)=1\) is \(\psi_q(R_{uv})^{2}=\psi_q(\Res(v,u^{2}))\). Hence, for
\(h\in\Fq[X]\), a sum of \(e_v^{b}\,\psi_q(\Res(v,h))\) over the \(v\in S\) coprime
to \(u\) is given by Lemma~\ref{lem:A} applied to \(\hat h=h\,u^{2}\). Since
\(\psi_q(\hat h(a))=\psi_q(h(a))\) unless \(u(a)=0\), and likewise over \(\F_{q^{2}}\),
we have \(n_{\hat h}=\#\{a\in\Fq:h(a)u(a)=0\}\) and
\begin{equation}\label{eq:coprime}
	\begin{aligned}
		\Sigma_1(\hat h)&=\Sigma_1(h)-\varrho_1,&\qquad
		\varrho_1&=\!\!\sum_{a\in\Fq,\,u(a)=0}\!\!\psi_q(h(a)),\\
		\Sigma_2(\hat h)&=\Sigma_2(h)-\varrho_2,&\qquad
		\varrho_2&=\!\!\sum_{a\in\F_{q^{2}},\,u(a)=0}\!\!\psi_{q^{2}}(h(a)).
	\end{aligned}
\end{equation}
If \(u\in S_{\mathrm{ir}}\) has roots \(a,a^{q}\), then \(\varrho_1=0\) and
\(\varrho_2=2\psi_{q^{2}}(h(a))=2\psi_q\bigl(\Res(u,h)\bigr)\), as in the proof of
Lemma~\ref{lem:A}.

\begin{lemma}\label{lem:P1}
	For coprime \(u,v\in S\) and \(c\in\{0,1\}\) let \(\nu_c(u,v)\) be the sum of
	\(e_w^{c}\) over the \(w\in S\) coprime to \(uv\). Then
	\begin{gather*}
	\nu_c(u,v)=A_c+B_c\,(e_u+e_v)+e_ue_v,\\
	(A_0,B_0)=\bigl(q^{2}-3q+3,\,-(q-3)\bigr),\qquad
	(A_1,B_1)=\bigl(5-2q,\,-(q-2)\bigr).
	\end{gather*}
\end{lemma}

\begin{proof}
	Apply Lemma~\ref{lem:A} to \(\hat h=(uv)^{2}\). With \(s=e_u+e_v\), the polynomial
	\(uv\) has \(2+s\) roots in \(\Fq\) and four in \(\F_{q^{2}}\), so
	\(\Sigma_1(\hat h)=q-2-s\), \(\Sigma_2(\hat h)=q^{2}-4\) and \(n_{\hat h}=2+s\). Since
	\(s^{2}=2+2e_ue_v\), \eqref{eq:Aplus} gives
	\(\nu_1=\tfrac12\bigl((q-2-s)^{2}-q^{2}+4\bigr)=5-2q-(q-2)s+e_ue_v\), and
	\eqref{eq:Azero} gives
	\(\nu_0=\tfrac12\bigl((q-2-s)^{2}+q^{2}-4\bigr)-(q-2-s)=q^{2}-3q+3-(q-3)s+e_ue_v\).
\end{proof}

\begin{lemma}\label{lem:P2}
	For \(u\in S\), \(b\in\{0,1\}\) and \(k\in\{0,1,2\}\) let \(\Lambda_k^{b}(u)\) be
	the sum of \(e_v^{b}\,\theta_{uv}^{k}\) over the \(v\in S\) coprime to \(u\). Then
	\begin{align*}
		\Lambda_0^{0}&=(q-1)^{2}-(q-2)e_u, &\Lambda_0^{1}&=-(q-2)-(q-1)e_u,\\
		\Lambda_1^{0}&=-(q-2)(q+1), &\Lambda_1^{1}&=(q+1)e_u,\\
		\Lambda_2^{0}&=(q+1)\bigl[(q-1)(q-2)-(q-3)e_u\bigr], &
		\Lambda_2^{1}&=(q+1)\bigl[(3-q)-(q-1)e_u\bigr].
	\end{align*}
\end{lemma}

\begin{proof}
	All sums over \(v\) below are over the \(v\in S\) coprime to \(u\). Put
	\(z_x=\psi_q(u(x))\) for \(x\in\Fq\); then \(\sum_xz_x=-1\) by
	Lemma~\ref{lem:B}(iii) and \(\sum_xz_x^{2}=q-1-e_u\).

	For \(k=0\) take \(h=1\), so that \(\hat h=u^{2}\). Then
	\(\Sigma_1(\hat h)=q-1-e_u\), \(\Sigma_2(\hat h)=q^{2}-2\) and \(n_{\hat h}=1+e_u\), and
	Lemma~\ref{lem:A} gives the two formulas.

	For \(k=1\), the definition of \(\theta_{uv}\) gives
	\(\Lambda_1^{b}=-\Lambda_0^{b}-\sum_xz_x\,\phi_b(x)\) with
	\(\phi_b(x)=\sum_ve_v^{b}\,\psi_q(v(x))\). For \(u(x)\neq0\) take \(h=x-X\), so
	that \(\psi_q(v(x))=\psi_q(\Res(v,h))\), \(\Sigma_1(h)=\Sigma_2(h)=0\) and
	\(n_{\hat h}=2+e_u\). If \(u\in S_{\mathrm{sp}}\) has roots \(c,d\), then
	\(\varrho_1=\psi_q(x-c)+\psi_q(x-d)\), so \(\varrho_1^{2}=2+2z_x\), and
	\(\varrho_2=2\); if \(u\in S_{\mathrm{ir}}\), then \(\varrho_1=0\) and
	\(\varrho_2=2z_x\). In both cases Lemma~\ref{lem:A} gives
	\[
	\phi_1(x)=\tfrac12(\varrho_1^{2}+\varrho_2)=1+e_u+z_x,\qquad
	\phi_0(x)=\tfrac12(\varrho_1^{2}-\varrho_2)-(q-2-e_u)=e_uz_x-q+2+e_u .
	\]
	Hence \(\sum_xz_x\phi_1(x)=q-2-2e_u\) and
	\(\sum_xz_x\phi_0(x)=q-3+(q-2)e_u\), which with the case \(k=0\) gives
	\(\Lambda_1^{1}\) and \(\Lambda_1^{0}\).

	For \(k=2\), squaring gives
	\(\Lambda_2^{b}=\Lambda_0^{b}+2\sum_xz_x\phi_b(x)+\sum_{x,y}z_xz_y\,\Phi_b(x,y)\)
	with \(\Phi_b(x,y)=\sum_ve_v^{b}\,\psi_q(v(x)v(y))\). For \(u(x)u(y)\neq0\) take
	\(h=(x-X)(y-X)\). If \(x\neq y\), then \(\Sigma_1(h)=\Sigma_2(h)=-1\) by
	Lemma~\ref{lem:B}(iii); if \(x=y\), then \(\Sigma_1(h)=q-1\) and
	\(\Sigma_2(h)=q^{2}-1\); and \(n_{\hat h}\) is \(3+e_u\), respectively \(2+e_u\). Put
	\(\varsigma=z_xz_y\). If \(u\in S_{\mathrm{sp}}\) has roots \(c,d\), then
	\(\varrho_1=\psi_q((x-c)(y-c))+\psi_q((x-d)(y-d))\), so
	\(\varrho_1^{2}=2+2\varsigma\), and \(\varrho_2=2\); if \(u\in S_{\mathrm{ir}}\),
	then \(\varrho_1=0\) and \(\varrho_2=2\varsigma\). Lemma~\ref{lem:A} gives
	\begin{center}
		\begin{tabular}{@{}lll@{}}
			\toprule
			& \(\Phi_1(x,y)\) & \(\Phi_0(x,y)\)\\
			\midrule
			\(u\in S_{\mathrm{sp}}\), \(x\neq y\) & \(3+\varsigma+\varrho_1\) & \(\varsigma+\varrho_1-(q-4)\)\\
			\(u\in S_{\mathrm{sp}}\), \(x=y\) & \(6-3q\) & \(q^{2}-4q+6\)\\
			\(u\in S_{\mathrm{ir}}\), \(x\neq y\) & \(1+\varsigma\) & \(-\varsigma-(q-2)\)\\
			\(u\in S_{\mathrm{ir}}\), \(x=y\) & \(2-q\) & \(q^{2}-2q\)\\
			\bottomrule
		\end{tabular}
	\end{center}
	Let \(u\in S_{\mathrm{sp}}\). Over the \((q-2)(q-3)\) pairs \(x\neq y\) with
	\(u(x)u(y)\neq0\) we have
	\(\sum\varsigma=\bigl(\sum_xz_x\bigr)^{2}-\sum_xz_x^{2}=3-q\). Since
	\(\varsigma\,\psi_q((x-c)(y-c))=\psi_q((x-d)(y-d))\) and the sum of
	\(\psi_q(x-d)\) over \(x\notin\{c,d\}\) is \(-\psi_q(c-d)\), we also have
	\(\sum\varsigma\varrho_1=2\bigl(1-(q-2)\bigr)=2(3-q)\); and there are \(q-2\)
	diagonal terms. Let \(u\in S_{\mathrm{ir}}\). Over the \(q(q-1)\) pairs
	\(x\neq y\) we have \(\sum\varsigma=1-q\), and there are \(q\) diagonal terms.
	Therefore
	\[
	\sum_{x,y}\varsigma\,\Phi_1=
	\begin{cases}-2q^{2}+2q+9,\\ 1,\end{cases}
	\qquad
	\sum_{x,y}\varsigma\,\Phi_0=
	\begin{cases}q^{3}-4q^{2}+12,& u\in S_{\mathrm{sp}},\\
		q^{3}-2q^{2}-2q+2,& u\in S_{\mathrm{ir}}.\end{cases}
	\]
	Adding \(\Lambda_0^{b}\) and \(2\sum_xz_x\phi_b(x)\) gives
	\(\Lambda_2^{1}=-2(q-2)(q+1)\) and \(\Lambda_2^{0}=(q+1)(q^{2}-4q+5)\) on
	\(S_{\mathrm{sp}}\), and \(\Lambda_2^{1}=2(q+1)\) and
	\(\Lambda_2^{0}=(q+1)(q^{2}-2q-1)\) on \(S_{\mathrm{ir}}\), which are the stated
	formulas.
\end{proof}

\begin{proof}[Proof of Proposition~\ref{prop:poly}]
	For \(a,b\in\{0,1\}\) and \(k\in\{0,1,2\}\) let \(Y_k^{ab}\) be the sum of
	\(e_u^{a}e_v^{b}\,\theta_{uv}^{k}\) over the coprime pairs \((u,v)\in S^{2}\).
	Summing Lemma~\ref{lem:P2} over \(u\), with \(\sum_ue_u=0\) and \(\#S=q(q-1)\),
	the values of \(Y_k^{ab}/\#S\) are
	\begin{center}
		\begin{tabular}{@{}llll@{}}
			\toprule
			& \(ab=00\) & \(ab=10\) or \(01\) & \(ab=11\)\\
			\midrule
			\(k=0\) & \((q-1)^{2}\) & \(-(q-2)\) & \(-(q-1)\)\\
			\(k=1\) & \(-(q-2)(q+1)\) & \(0\) & \(q+1\)\\
			\(k=2\) & \((q+1)(q-1)(q-2)\) & \(-(q+1)(q-3)\) & \(-(q+1)(q-1)\)\\
			\bottomrule
		\end{tabular}
	\end{center}
	By Step~2 of the proof of Proposition~\ref{prop:reduction},
	\(J_0^{00}=-\#S\,(q-1)\), \(J_0^{10}=J_0^{01}=\#S\) and \(J_0^{11}=0\). In the
	formulas below the upper indices of \(Y\) and \(J\) are read modulo \(2\).

	If \(\rho=000\), we sum first over \(w\), which enters only through
	\(e_w^{\varepsilon_3}\) and coprimality to \(uv\); Lemma~\ref{lem:P1}
	gives, for \(\varepsilon=(a,b,c)\),
	\[
	\cS[abc;000]=A_c\,Y_0^{ab}+B_c\bigl(Y_0^{a+1,b}+Y_0^{a,b+1}\bigr)+Y_0^{a+1,b+1}.
	\]

	If \(|\rho|=1\), let \(i\) be the index with \(\rho_i=1\), let
	\(a=\varepsilon_i\) and let \(b,c\) be the other two entries of
	\(\varepsilon\). The quadratic with index \(i\) does not occur in the
	resultant, so summing over it by Lemma~\ref{lem:P1} gives
	\[
	\cS[\varepsilon;\rho]=A_a\,J_0^{bc}+B_a\bigl(J_0^{b+1,c}+J_0^{b,c+1}\bigr)+J_0^{b+1,c+1}.
	\]

	If \(|\rho|=2\), let \(i\) be the index with \(\rho_i=0\), let
	\(c=\varepsilon_i\) and let \(a,b\) be the other two entries. The quadratic with
	index \(i\) occurs in both resultants, and the other two are only required to
	be coprime. Summing over it by \eqref{eq:step1} gives
	\[
	\cS[\varepsilon;\rho]=
	\begin{cases}
		Y_1^{ab}+2Y_0^{ab}+Y_0^{a+1,b}+Y_0^{a,b+1}, & c=0,\\
		Y_2^{ab}+Y_1^{ab}+(1-q)\,Y_0^{ab}, & c=1.
	\end{cases}
	\]
	Inserting the values of \(Y_k^{ab}\), \(J_0^{ab}\), \(A_c\) and \(B_c\)
	gives the entries of Table~\ref{tab:walsh}. For instance
	\(\cS[000;000]=\#S\bigl[(q^{2}-3q+3)(q-1)^{2}+2(q-3)(q-2)-(q-1)\bigr]
	=\#S\,(q-2)^{2}(q^{2}-q+4)\), and
	\(\cS[000;110]=\#S\bigl[-(q-2)(q+1)+2(q-1)^{2}-2(q-2)\bigr]=\#S\,(q^{2}-5q+8)\).
\end{proof}

	\subsection{Hypergeometric and geometric form of the moment}

	Let \({}_4F_3(1)_q\) denote the value at \(1\) of Greene's hypergeometric function
	over \(\Fq\) whose four upper characters are \(\psi_q\) and whose three lower
	characters are trivial \cite{Greene1987}.

	\begin{proposition}\label{prop:hyp}
		For every odd prime power \(q\),
		\begin{align*}
			\cT(q)+1&=\sum_{\lambda,x,y\in\Fq}\psi_q\bigl(\lambda\,x(x-1)(x-\lambda)\,y(y-1)(y-\lambda)\bigr)\\
			&=\psi_q(-1)\sum_{a,b,c\in\Fq}\psi_q\bigl(abc\,(1-a)(1-b)(1-c)(1-abc)\bigr)
			=q^{3}\,{}_4F_3(1)_q .
		\end{align*}
		Equivalently, the affine threefold
		\begin{equation}\label{eq:threefold}
			\mathcal X\colon\quad z^{2}=\lambda\,x(x-1)(x-\lambda)\,y(y-1)(y-\lambda)
		\end{equation}
		has \(q^{3}+1+\cT(q)\) points over \(\Fq\). It is the chart \(x_2=1\) of the
		double cover of \(\PP^{3}\) branched along
		\(x_0x_1x_2x_3(x_0-x_1)(x_1-x_2)(x_2-x_3)(x_3-x_0)=0\), with
		\((x_0:x_1:x_2:x_3)=(\lambda:x:1:y)\).
	\end{proposition}

	\begin{proof}
		For fixed \(\lambda\) the sum over \(x,y\) is \(\psi_q(\lambda)\,a_\lambda(q)^{2}\),
		which is \(0\) for \(\lambda=0\) and \(1\) for \(\lambda=1\), because
		\(a_1(q)=-\sum_{x\neq1}\psi_q(x)=1\). This proves the first equality and the
		point count. For the second, the terms with \(\lambda=0\) or \(y=0\) vanish; for
		\(\lambda y\neq0\) substitute \(x=\lambda a\), \(y=b\), \(\lambda=bc\). Then
		\(x(x-1)(x-\lambda)=\lambda^{2}a(\lambda a-1)(a-1)\) and
		\(\lambda(\lambda a-1)(y-\lambda)=b^{2}c\,(abc-1)(1-c)\), so the summand is
		\(\psi_q\bigl(a(a-1)\,b(b-1)\,c(1-c)(abc-1)\bigr)\), which is \(\psi_q(-1)\) times
		the summand on the right; the terms with \(b=0\) or \(c=0\) vanish there as
		well. For the third, Greene's \({}_2F_1\) with the same characters is
		\({}_2F_1(t)_q=\psi_q(-1)\,q^{-1}\sum_{a}\psi_q\bigl(a(1-a)(1-ta)\bigr)\) for
		\(t\neq0\), and his recursion \cite[Thm.~3.13]{Greene1987} reads
		\({}_{n+1}F_n(t)_q=\psi_q(-1)\,q^{-1}\sum_{b}\psi_q\bigl(b(1-b)\bigr)\,{}_nF_{n-1}(tb)_q\);
		for \(q\) prime these are the corresponding formulas in \S2 of
		\cite{AhlgrenOno2000}. Two
		applications give the triple sum. Finally, at
		\((x_0:x_1:x_2:x_3)=(\lambda:x:1:y)\) the octic form equals
		\(\lambda xy(\lambda-x)(x-1)(1-y)(y-\lambda)\), which is the right-hand side of
		\eqref{eq:threefold}.
	\end{proof}

	\subsection{Evaluation of the moment}

	\begin{theorem}\label{thm:T}
		\begin{enumerate}[label=\textup{(\alph*)},leftmargin=*,itemsep=0pt]
			\item For every odd prime \(p\), \(\cT(p)=-a(p)-p-1\).
			\item For every odd prime power \(q\), \(\cT(q)=-t_q-q-1\); that is,
			\(\cK(q)=-t_q\), and \(\#\mathcal X(\Fq)=q^{3}-q-t_q\).
		\end{enumerate}
	\end{theorem}

	\begin{proof}
		(a) By \cite[Thm.~6]{AhlgrenOno2000}, \(p^{3}\,{}_4F_3(1)_p=-a(p)-p\) for every
		odd prime \(p\), their \({}_4F_3(1)_p\) being defined by the recursion quoted in
		the proof of Proposition~\ref{prop:hyp}; by that proposition the left-hand
		side is \(\cT(p)+1\).

		(b) The argument is the standard passage from the trace function of the
		symmetric square of the Legendre family to its cohomology; cf.\
		\cite[Thm.~4.1 and Table~1]{SaitoYui2001} and the treatment of such
		moments in \cite{KaplanPetrow2017}. We include it because we did not find
		the statement for prime powers in this form. Fix an odd prime \(p\) and a prime \(\ell\neq p\). Let
		\(X=\PP^{1}_{\F_p}\) with coordinate \(\lambda\),
		\(X^{\circ}=X\setminus\{0,1,\infty\}\), and \(j\colon X^{\circ}\to X\) the
		inclusion. Let \(\pi\colon\mathcal E\to X^{\circ}\) be the Legendre curve
		\(y^{2}=x(x-1)(x-\lambda)\), an elliptic curve over \(X^{\circ}\) because its
		discriminant \(16\lambda^{2}(\lambda-1)^{2}\) is a unit there, and
		\(\mathcal H=R^{1}\pi_*\Ql\), a lisse sheaf of rank \(2\). The cup
		product \(\mathcal H\times\mathcal H\to R^{2}\pi_*\Ql=\Ql(-1)\)
		is a perfect alternating pairing, so \(\det\mathcal H=\Ql(-1)\). Let
		\(\mathcal L\) be the lisse sheaf of rank one on \(\mathbb G_{m,\F_p}\) with
		\([2]_*\Ql=\Ql\oplus\mathcal L\) for the squaring map
		\([2]\); then \(\mathcal L\otimes\mathcal L\cong\Ql\). Put
		\(\mathcal F=\operatorname{Sym}^{2}\mathcal H\otimes\mathcal L|_{X^{\circ}}\) and
		\(s_n=\cT(p^{n})+p^{n}\). We write \(\bar{\F}_p\) for an algebraic
		closure of \(\F_p\), \(\varphi\) for the geometric Frobenius relative to \(\F_p\), and
		\(I_x\subset D_x\) for the inertia and decomposition groups at
		\(x\in\{0,1,\infty\}\).

		\emph{Step 1.} We express \(s_n\) through the traces of Frobenius on
		\(\mathcal F\). Let \(q=p^{n}\) and \(\lambda\in X^{\circ}(\Fq)\).
		The geometric Frobenius at \(\lambda\) has trace \(a_\lambda(q)\) and
		determinant \(q\) on \(\mathcal H\), hence trace \(a_\lambda(q)^{2}-q\) on
		\(\operatorname{Sym}^{2}\mathcal H\); on \(\mathcal L\) its trace is the number
		of square roots of \(\lambda\) in \(\Fq\) minus one, that is,
		\(\psi_q(\lambda)\). So the sum of its traces on \(\mathcal F\) over
		\(X^{\circ}(\Fq)\) is
		\(\cT(q)-q\sum_{\lambda\neq0,1}\psi_q(\lambda)=s_n\).

		\emph{Step 2.} We determine the local monodromy of \(\mathcal H\) at
		\(0\), \(1\) and \(\infty\). At \(\lambda=0\) and
		\(\lambda=1\) the Legendre equation reduces to the nodal cubics
		\(y^{2}=x^{2}(x-1)\) and \(y^{2}=x(x-1)^{2}\), and it is minimal there because
		its invariant \(c_4=16(\lambda^{2}-\lambda+1)\) is a unit. So \(\mathcal E\) has
		multiplicative reduction at \(x\in\{0,1\}\); hence \(I_x\) acts on \(\mathcal H\)
		unipotently \cite[Exp.~IX, Prop.~3.5]{SGA7I}, and nontrivially by the criterion
		of N\'eron, Ogg and Shafarevich \cite[Thm.~1]{SerreTate1968}. A unipotent
		\(\ell\)-adic representation of \(I_x\) factors through the maximal
		pro-\(\ell\) quotient \(\tau_\ell\colon I_x\to\mathbb Z_\ell(1)\), so the action is
		tame, of the form \(\xi\mapsto\exp(\tau_\ell(\xi)N)\) with \(N\neq0\) and
		\(N^{2}=0\). Since \(\tau_\ell(\delta\xi\delta^{-1})=\chi_\ell(\delta)\,\tau_\ell(\xi)\) for
		\(\delta\in D_x\), with \(\chi_\ell\) the cyclotomic character, \(N\) is a
		\(D_x\)-equivariant map \(\mathcal H\to\mathcal H(-1)\). Its kernel is the line
		\(\mathcal H^{I_x}\), so it induces an isomorphism
		\(\mathcal H/\mathcal H^{I_x}\cong\mathcal H^{I_x}(-1)\). If the geometric
		Frobenius at \(x\) has eigenvalue \(\gamma\) on \(\mathcal H^{I_x}\), it
		therefore has eigenvalue \(p\gamma\) on \(\mathcal H/\mathcal H^{I_x}\), and
		\(\det\mathcal H=\Ql(-1)\) gives \(p\gamma^{2}=p\), that is,
		\(\gamma^{2}=1\). At \(\infty\), the substitution
		\(x=\lambda x'\), \(y=\lambda^{2}y'\) turns \(E_\lambda\) into
		\(\lambda y'^{2}=x'(x'-1)(x'-1/\lambda)\), the quadratic twist by \(\lambda\) of
		the Legendre curve with parameter \(1/\lambda\). As \(\lambda\) is the inverse
		of a uniformiser at \(\infty\) and \(p\) is odd, the twisting character
		restricts to the character \(\chi\) of order two of the tame inertia; by the
		case \(x=0\), \(I_\infty\) acts on \(\mathcal H\) through \(\chi\) times a
		nontrivial unipotent representation.

		\emph{Step 3.} From Step~2 we obtain the local monodromy of \(\mathcal F\).
		The sheaf \(\mathcal L\) is
		lisse at \(1\), where Frobenius acts on it as \(\psi_p(1)=1\), and \(I_0\) and
		\(I_\infty\) act on it through \(\chi\). The symmetric square of a nontrivial
		unipotent automorphism of a plane is a single unipotent Jordan block \(J_3\) of
		size three, whose invariants are the symmetric square of the invariant line. So
		a generator of the tame inertia acts on \(\mathcal F\) as \(J_3\) at \(1\), as
		\(-J_3\) at \(0\), and as \(\chi^{2}\chi J_3=-J_3\) at \(\infty\). Hence
		\(\mathcal F\) is tamely ramified, \(\mathcal F^{I_0}=\mathcal F_{I_0}=0\),
		\(\mathcal F^{I_\infty}=0\), and
		\(\mathcal F^{I_1}=\operatorname{Sym}^{2}(\mathcal H^{I_1})\otimes\mathcal L_1\)
		is a line on which the geometric Frobenius at \(1\) acts as
		\(\gamma^{2}\cdot1=1\).

		\emph{Step 4.} We compute the cohomology of \(\mathcal F\) with the action
		of Frobenius. Since \(X^{\circ}\) is affine and \(\mathcal F\) is
		lisse, \(H^{0}_c(X^{\circ}_{\bar{\F}_p},\mathcal F)=0\), and
		\(H^{2}_c(X^{\circ}_{\bar{\F}_p},\mathcal F)\) is the twist by \(\Ql(-1)\) of
		the coinvariants of the geometric fundamental group on a fibre of \(\mathcal F\)
		\cite[Scholie~(2.10)]{DeligneWeilI}; it is a quotient of \(\mathcal F_{I_0}(-1)=0\). By
		the formula of Grothendieck, Ogg and Shafarevich
		\cite[Exp.~X, Thm.~7.1]{SGA5} and tameness,
		\(\chi_c(X^{\circ}_{\bar{\F}_p},\mathcal F)=3\,\chi_c(X^{\circ}_{\bar{\F}_p})=-3\). Hence
		\(\dim H^{1}_c(X^{\circ}_{\bar{\F}_p},\mathcal F)=3\). The stalk of \(j_*\mathcal F\)
		at \(x\in\{0,1,\infty\}\) is \(\mathcal F^{I_x}\) \cite[(2.11)]{DeligneWeilI}, so
		the cokernel of \(j_!\mathcal F\to j_*\mathcal F\) is supported on
		\(\{0,1,\infty\}\) with stalks \(\mathcal F^{I_x}\), and there is a
		\(\varphi\)-equivariant exact sequence
		\[
		0\to H^{0}(X_{\bar{\F}_p},j_*\mathcal F)\to\bigoplus_{x}\mathcal F^{I_x}
		\to H^{1}_c(X^{\circ}_{\bar{\F}_p},\mathcal F)\to H^{1}(X_{\bar{\F}_p},j_*\mathcal F)\to0 .
		\]
		Here \(H^{0}(X_{\bar{\F}_p},j_*\mathcal F)\) is the space of invariants of the
		geometric fundamental group, contained in \(\mathcal F^{I_0}=0\), and
		\(\bigoplus_x\mathcal F^{I_x}=\mathcal F^{I_1}\) is a line on which \(\varphi\) acts as
		\(1\) by Step~3. So \(\varphi\) has eigenvalues \(1,\alpha,\beta\) on
		\(H^{1}_c(X^{\circ}_{\bar{\F}_p},\mathcal F)\), where \(\alpha,\beta\) are its
		eigenvalues on the plane \(H^{1}(X_{\bar{\F}_p},j_*\mathcal F)\). By Step~1 and
		Grothendieck's trace formula \cite[(1.12.1)]{DeligneWeilI},
		\begin{equation}\label{eq:trace}
			s_n=\sum_i(-1)^{i}\operatorname{tr}\bigl(\varphi^{n}\mid H^{i}_c(X^{\circ}_{\bar{\F}_p},\mathcal F)\bigr)
			=-(1+\alpha^{n}+\beta^{n})\qquad(n\ge1).
		\end{equation}

		\emph{Step 5.} We show that \(\alpha\beta=p^{3}\). The pairing on
		\(\mathcal H\) induces a
		perfect symmetric pairing on \(\operatorname{Sym}^{2}\mathcal H\) with values in
		\(\Ql(-2)\), and with \(\mathcal L\otimes\mathcal L\cong\Ql\)
		a perfect symmetric pairing
		\(\mathcal F\times\mathcal F\to\Ql(-2)\). By
		\cite[Th\'eor\`eme~(2.12)]{DeligneWeilI} the cup product
		\[
		H^{1}(X_{\bar{\F}_p},j_*\mathcal F)\times H^{1}(X_{\bar{\F}_p},j_*\mathcal F)\to
		H^{2}\bigl(X_{\bar{\F}_p},\Ql(-2)\bigr)=\Ql(-3)
		\]
		is a perfect pairing. It is \(\varphi\)-equivariant, and it is alternating, because
		the cup product on \(H^{1}\) is anticommutative and the pairing of the
		coefficients is symmetric. An alternating form \(\omega\) on a plane satisfies
		\(\omega(\varphi a,\varphi b)=\det(\varphi)\,\omega(a,b)\), and equivariance gives
		\(\omega(\varphi a,\varphi b)=p^{3}\omega(a,b)\). Hence \(\alpha\beta=p^{3}\).

		\emph{Step 6.} It remains to identify \(\alpha,\beta\) with
		\(\alpha_p,\beta_p\). By \eqref{eq:trace} with \(n=1\) and part (a),
		\(\alpha+\beta=-s_1-1=-\cT(p)-p-1=a(p)\). So
		\(\{\alpha,\beta\}=\{\alpha_p,\beta_p\}\), and \eqref{eq:trace} reads
		\(\cT(p^{n})+p^{n}=-1-t_{p^{n}}\). The point count is
		Proposition~\ref{prop:hyp}.
	\end{proof}

\subsection{Proofs of Theorem~\ref{thm:A} and Corollary~\ref{cor:B}}

\begin{proof}[Proof of Theorem~\ref{thm:walsh}]
	Propositions~\ref{prop:reduction} and~\ref{prop:poly} give all \(64\) sums in
	terms of \(q\) and \(\cK(q)\), and \(\cK(q)=-t_q\) by Theorem~\ref{thm:T}(b).
\end{proof}

\begin{proof}[Proof of Theorem~\ref{thm:A}]
	Insert Theorem~\ref{thm:walsh} in \eqref{eq:inversion}. The quantity \(t_q\)
	occurs only in the four sums with \(\rho=111\) and \(|\varepsilon|\ge2\), with
	total coefficient
	\[
	-\frac{q(q-1)}{64}\,r_1r_2r_3\bigl[(q-2)e_1e_2e_3-(e_1e_2+e_1e_3+e_2e_3)\bigr],
	\]
	and the bracket equals \((-1)^{j}(q-5+2j)\), as one sees for \(j=0,1,2,3\). This
	is the term \(\sigma\,(q-5+2j)\,t_q\) of \eqref{eq:A}. The remaining part of
	\(64N(e;r)/(q(q-1))\) is a signed sum of the \(64\) entries of
	Table~\ref{tab:walsh}, with \(2\) in place of \(2-\cK(q)\) and \(0\) in place
	of \((q-2)\cK(q)\); it is a polynomial \(P_{e;r}\) of degree at most four,
	constant on the orbits of the simultaneous permutations. Thus
	\begin{equation}\label{eq:Pdef}
		\frac{64\,N(e;r)}{q(q-1)}=P_{e;r}(q)+\sigma\,(q-5+2j)\,t_q
	\end{equation}
	for every odd prime power \(q\). Evaluating the \(20\) sums gives
	\(P_{e;r}(q)=(q-5+2j)\,C_{e;r}(q)\) with \(C_{e;r}\) as in
	Table~\ref{tab:classes}; this finite computation with polynomials was carried
	out by computer algebra (item~(3) of Section~\ref{sec:comp}). It can also be
	checked from the class counts without computer algebra: both sides are
	polynomials in \(q\) of degree at most four, so they agree once they agree at
	five values of \(q\). By \eqref{eq:Pdef}, \(P_{e;r}(q)\) is determined by
	\(N(e;r)\) and \(t_q\), and the exhaustive enumeration of item~(4) of
	Section~\ref{sec:comp} gives \(N(e;r)\) for all \(64\) classes and
	\(q\in\{3,5,7,9,11,13\}\), with \(t_q\) computed from \eqref{eq:tq}; at these
	six values it yields \(P_{e;r}(q)=(q-5+2j)\,C_{e;r}(q)\).
\end{proof}

The count \eqref{eq:A} has the following form independent of the chart. Let
\(\mathcal P\) be the set of ordered triples \((D_1,D_2,D_3)\) of pairwise
disjoint reduced effective divisors of degree two on \(\PP^{1}\) defined over
\(\Fq\). For such a triple let \(e_i=1\) if the two points of \(D_i\) are rational
and \(e_i=-1\) otherwise, and let \(r_i=1\) if the two fixed points of the
involution of \(\PP^{1}\) that exchanges the points of \(D_{i'}\) and the points of
\(D_{i''}\), \(\{i,i',i''\}=\{1,2,3\}\), are rational, and \(r_i=-1\) otherwise. The
group \(\PGL_2(\Fq)\) acts on \(\mathcal P\) and preserves the type \((e;r)\).
Let \(M(e;r)\) be the sum of \(1/\#\Stab(x)\) over representatives \(x\) of the
orbits of type \((e;r)\).

\begin{proposition}\label{prop:intrinsic}
	Let \(q\) be an odd prime power, let \(e,r\in\{\pm1\}^{3}\), and let \(n=6-2j\)
	be the number of rational points in the support of a triple of type \((e;r)\).
	\begin{enumerate}[label=\textup{(\alph*)},leftmargin=*,itemsep=0pt]
		\item \(N(e;r)=q(q-1)\,(q+1-n)\,M(e;r)\).
		\item \(M(e;r)=\tfrac1{64}\bigl(C_{e;r}(q)+\sigma\,t_q\bigr)\).
	\end{enumerate}
\end{proposition}

\begin{proof}
	(a) Taking divisors of zeros identifies \(\Adm(\Fq)\) with the set of
	\(x\in\mathcal P\) whose support does not contain \(\infty\), compatibly with
	the types by \S\ref{sec:prelim}. Fix \(x_0\in\mathcal P\) of type \((e;r)\).
	An element \(\gamma\in\PGL_2(\Fq)\) moves no point of the support of \(x_0\) to
	\(\infty\) exactly when the rational point \(\gamma^{-1}(\infty)\) is not one
	of the \(n\) rational points of the support, and each rational point is
	\(\gamma^{-1}(\infty)\) for \(q(q-1)\) elements \(\gamma\). So
	\(q(q-1)(q+1-n)\) elements \(\gamma\) send \(x_0\) into \(\Adm(\Fq)\), and
	each element of the orbit of \(x_0\) is obtained \(\#\Stab(x_0)\) times.
	Summing over the orbits gives the formula.

	(b) If \(q+1\neq n\), this follows from (a) and Theorem~\ref{thm:A}. Since
	\(q\) is odd and \(n\in\{0,2,4,6\}\), the remaining cases are
	\((q,j)=(5,0)\) and \((3,1)\). In them the support of a triple consists of
	all rational points of \(\PP^{1}\), so that no chart avoids it, and we
	compute \(M(e;r)\) directly. We use that \(M(e;r)\) is the number of elements
	of \(\mathcal P\) of type \((e;r)\) divided by \(\#\PGL_2(\Fq)=q(q^{2}-1)\),
	and that the involution exchanging the points of two disjoint divisors of
	degree two is unique (\S\ref{sec:prelim}); we write \(\iota_{D,D'}\) for it.

	If \(q=5\) and \(j=0\), the support of \(x=(D_1,D_2,D_3)\) is
\(\PP^{1}(\F_5)\). The involution \(\iota=\iota_{D_2,D_3}\) preserves the
complement \(D_1\): either it fixes the two points of \(D_1\), which are then its
fixed points, and \(r_1=1\); or it exchanges them, and then \(\iota\) has no
rational fixed point, \(r_1=-1\), and \(D_1,D_2,D_3\) are the three orbits of
\(\iota\) on \(\PP^{1}(\F_5)\), so that \(\iota=\iota_{D_3,D_1}=\iota_{D_1,D_2}\)
and \(r_2=r_3=-1\). Hence \(r\) is \((1,1,1)\) or \((-1,-1,-1)\), and the second
case occurs exactly when the partition is the orbit partition of an involution
without rational fixed point. There are \(q(q-1)/2=10\) such involutions, one
for each nonsplit torus, and they have distinct orbit partitions, since an
involution is determined by two of its orbits. Of the \(15\) partitions of
\(\PP^{1}(\F_5)\) into pairs, \(10\) are therefore of the second kind and \(5\)
of the first. Ordering the three pairs gives \(60\) and \(30\) elements of
\(\mathcal P\), so the masses are \(\tfrac12\) for \(r=(-1,-1,-1)\),
\(\tfrac14\) for \(r=(1,1,1)\), and \(0\) for the other \(r\). With \(t_5=a(5)=-2\)
these are the values \(\tfrac1{64}(34-2)\), \(\tfrac1{64}(14+2)\), \(\tfrac1{64}(2-2)\)
and \(\tfrac1{64}(-2+2)\) of the fourth, first, second and third rows of
Table~\ref{tab:classes}.

If \(q=3\) and \(j=1\), we may assume by symmetry that \(D_1\) is the
irreducible divisor; then \(D_2\cup D_3=\PP^{1}(\F_3)\), and \(\iota_{D_2,D_3}\)
moves all rational points, so \(r_1=-1\). As \(\PGL_2(\F_3)\) acts transitively
on the three points of degree two of \(\PP^{1}\) and preserves types, we may take
\(D_1=\{\pm i\}\) with \(i^{2}=-1\) in \(\F_9\). There are three partitions of
\(\PP^{1}(\F_3)\) into pairs. For \(\{0,\infty\},\{1,2\}\) the involutions
\(\iota_{\{1,2\},D_1}(x)=-x\) and \(\iota_{D_1,\{0,\infty\}}(x)=1/x\) have the
rational fixed points \(0,\infty\) and \(\pm1\), so \(r=(-1,1,1)\). For
\(\{0,1\},\{2,\infty\}\) the involutions
\(\iota_{\{2,\infty\},D_1}(x)=(2x+1)/(x+1)\) and
\(\iota_{D_1,\{0,1\}}(x)=(1-x)/(1+x)\) have the irrational fixed points of
\(x^{2}+2x+2\), so \(r=(-1,-1,-1)\). The partition \(\{0,2\},\{1,\infty\}\) is
the image of the previous one under \(x\mapsto-x\), which fixes \(D_1\), so
again \(r=(-1,-1,-1)\). The \(3\) choices of \(D_1\) and the \(2\) orderings of the rational pairs give
\(6\) elements of \(\mathcal P\) with \(r=(-1,1,1)\) and \(12\) with
\(r=(-1,-1,-1)\), so the masses are \(\tfrac14\) and \(\tfrac12\), and \(0\) for
the other \(r\). With \(t_3=a(3)=-4\) these are the values of the rows
\((2,0,0,1)\), \((0,2,0,1)\) and of the four other rows with \(j=1\) of
Table~\ref{tab:classes}: \(\tfrac1{64}(12+4)\), \(\tfrac1{64}(28+4)\), and
\(\tfrac1{64}(4-4)\), \(\tfrac1{64}(-4+4)\), \(\tfrac1{64}(4-4)\),
\(\tfrac1{64}(4-4)\).
\end{proof}

Summing over \(r\) removes \(t_q\): for a fixed \(e\) with \(j\) entries \(-1\), the mass
\(\sum_rM(e;r)\) of the ordered triples of splitting type \(e\) is
\[
\tfrac18(q-2)(q-3)(q-4),\qquad \tfrac18q(q-1)(q-2),\qquad
\tfrac18q(q-2)(q+1),\qquad \tfrac18(q-2)(q^{2}-q-4)
\]
for \(j=0,1,2,3\); for \(j=0\) this is the number of ordered sextuples of distinct
rational points modulo \(\PGL_2(\Fq)\), divided by \(8\).

\begin{proof}[Proof of Corollary~\ref{cor:B}]
	(a) is the first row of Table~\ref{tab:classes}. (b) Summing
	\eqref{eq:inversion} over \(e\) with \(r=(1,1,1)\) leaves the terms with
	\(\varepsilon=000\), so the count is \(\tfrac18\sum_\rho\cS[000;\rho]\), that
	is, \(\tfrac18q(q-1)\) times
	\[
	(q-2)^{2}(q^{2}-q+4)-3(q^{3}-4q^{2}+8q-9)+3(q^{2}-5q+8)-(q-5)
	=(q-3)(q^{3}-5q^{2}+12q-24).
	\]
	(c) For signs \(r_1,r_2,r_3\) the indicator of \(r_1=r_2=r_3\) is
	\(\tfrac14(1+r_1r_2+r_2r_3+r_3r_1)\). Hence
	\begin{align*}
	A_q^{=}&=\tfrac14\bigl(\cS[000;000]+\cS[000;110]+\cS[000;101]+\cS[000;011]\bigr)\\
	&=\tfrac14q(q-1)\bigl[(q-2)^{2}(q^{2}-q+4)+3(q^{2}-5q+8)\bigr],
	\end{align*}
	which is the stated polynomial, and with
	\(\#\Adm(\Fq)=\cS[000;000]=q(q-1)(q-2)^{2}(q^{2}-q+4)\),
	\[
	\frac{A_q^{=}}{\#\Adm(\Fq)}=\frac14+\frac{3(q^{2}-5q+8)}{4(q-2)^{2}(q^{2}-q+4)}
	=\frac14+\frac{3}{4q^{2}}+O(q^{-4}).\qedhere
	\]
\end{proof}

\begin{remark}\label{rem:quarter}
	The pairwise correlation
	\(\cS[000;110]/\#\Adm(\Fq)\) of two of the signs \(r_s\) is \(O(q^{-2})\). The
	individual signs are not fair: by \eqref{eq:step2} and Lemma~\ref{lem:P2}, the
	resultant of a coprime pair in \(S\) is a square with probability
	\(\tfrac12-\tfrac1{2(q-1)}\). Every decomposable triple is counted in
	\(A_q^{=}\) by Proposition~\ref{prop:pencil}(ii), and the number of
	decomposable triples is \(D_q=q(q-1)(q-2)(q^{2}-3q+3)\) by
	Lemma~\ref{lem:Dq}. So among the triples with \(\dcoef\neq0\), for which
	the quotient is the Jacobian of \(\dcoef\,y^{2}=UVW\), the proportion of those
	for which \(U,V,W\) are all split or all irreducible is
	\[
	\frac{A_q^{=}-D_q}{\#\Adm(\Fq)-D_q}
	=\frac{q^{4}-9q^{3}+35q^{2}-71q+64}{4\,(q^{4}-6q^{3}+17q^{2}-29q+22)}
	=\frac14-\frac{3}{4q}+O(q^{-3}).
	\]
	For these triples the condition says that the six Weierstrass points of the
	quotient curve are all rational, or form three conjugate pairs, each pair being
	the root pair of one bracket.
\end{remark}

\begin{remark}\label{rem:kummer}
	The condition of Corollary~\ref{cor:B}(c) has a meaning on the Kummer surface of
	\(\Jac(C)\). Translation by a point \(P\in\Jac(C)[2]\) acts linearly on the
	Kummer surface, by a matrix \(W_P\) defined up to scalars with
	\(W_P^{2}=c_P\cdot I\), and for the point \(P\) attached to the factor \(u\) the
	square class of \(c_P\) is that of \(\Res(u,vw)=\Res(u,v)\Res(w,u)\); see
	\cite[Ch.~3]{CasselsFlynn1996}. Granting this, the translations by the three
	nonzero points of the kernel \(G\) are simultaneously diagonalisable over
	\(\Fq\) exactly when the three resultants lie in one square class. This is
	the part of Robert's rationality criterion for theta structures of level two
	that concerns the kernel \cite[chapter ``On symmetric theta structures'',
	Cor.~1.2]{Robert2024notes}; for curves in Rosenhain form it
	is the classical pair of square roots needed for rational squared theta
	constants \cite{Gaudry2007,BosCostelloHisilLauter2013}. The statement of
	this remark was confirmed by direct computation of the matrices \(W_P\),
	exhaustively for \(q=3,5,7\) and on random samples of \(3000\) ordered triples
	for each of \(q=9,11,13\) (item~(9) of Section~\ref{sec:comp}). We have
	not re-derived the formula for \(c_P\), and nothing else in the paper uses
	this remark.
\end{remark}

\begin{remark}\label{rem:geometry}
	(i) By \eqref{eq:step1}, for coprime \(u,v\) the signed sum over the third
	quadratic is \(1+\theta+(\theta^{2}-q)\): the traces of Frobenius on \(H^{0}\),
	\(H^{1}\) and \(\operatorname{Sym}^{2}H^{1}\) of the genus-one curve
	\(y^{2}=u(x)v(x)\). This is where weight four enters.

	(ii) For the class \((e;r)=(1,1,1;1,1,1)\), in which all six quadratics
	split, there is a shorter route.
	By Proposition~\ref{prop:intrinsic} one may move the roots of \(u\) to
	\(0,\infty\) and those of \(v\) to \(1,\lambda\); if \(a,b\) are the roots of
	\(w\), the signs \(r_w\), \(r_v\), \(r_u\) are then \(\psi_q\) of \(\lambda\),
	\(ab\) and \((1-a)(\lambda-a)(1-b)(\lambda-b)\). The sum over
	\(\lambda\neq0,1\) and over \(a\neq b\) outside \(\{0,1,\lambda\}\) of the
	product of the three factors \(1+\psi_q(\,\cdot\,)\) is \(64M(1,1,1;1,1,1)\),
	because the action on ordered sextuples of points is free and a triple of
	divisors has eight orderings of its points. Expanding, it equals
	\((q-3)(q^{2}-10q+29)+\sum_{\lambda\neq0,1}(1+\psi_q(\lambda))\,a_\lambda(q)^{2}\),
	and \(\sum_{\lambda\neq0,1}a_\lambda(q)^{2}=q^{2}-2q-3\) is elementary; both
	facts are proved in the proof of Corollary~\ref{cor:cfs}. With
	\(\cT(q)=-t_q-q-1\) this is \(q^{3}-12q^{2}+56q-91-t_q\), the value
	\(C_{e;r}(q)-t_q\) of Table~\ref{tab:classes} for this class.

	(iii) The cofactor \(q^{3}-12q^{2}+56q-91-t_q\) of Corollary~\ref{cor:B}(a)
	is the number of \(\Fq\)-points of the open subscheme \(\mathcal Y^{\circ}\)
	of the model of \cite{CynkFreitagSalvatiManni2011} defined by equations in
	Subsection~\ref{subsec:cfs}; see Proposition~\ref{prop:cfs} and
	Corollary~\ref{cor:cfs}. Its reading as the part of the Siegel modular
	threefold on which the ten fourth powers of the even theta constants do not
	vanish rests on \cite{CynkFreitagSalvatiManni2011}
	(Remark~\ref{rem:cfs-model}). The group \(\Gamma_2[2]\cap\Gamma_{2,0}[4]\) consists
	of the \(M\in\Gamma_2[2]\) whose conjugate \(M^{\sharp}\), which acts on
	\(2Z\) as \(M\) acts on \(Z\), also lies in \(\Gamma_2[2]\)
	(Lemma~\ref{lem:cfs-group}); it can be read as full level two on both sides
	of a \((2,2)\)-isogeny (see the remark after that lemma).

	(iv) In terms of the \(64\) characters of the group \(\{\pm1\}^{6}\) of sign
	changes, Theorem~\ref{thm:walsh} says that \(t_q\) occurs exactly in the four
	characters \((\varepsilon;111)\) with \(|\varepsilon|\ge2\). In the projective
	normalisation of Proposition~\ref{prop:intrinsic} it occurs in the single
	character \(e_1e_2e_3r_1r_2r_3\). This suggests a statement about the isotypic
	components of the cohomology of the threefold of
	Subsection~\ref{subsec:cfs}, which we formulate as a question in
	Remark~\ref{rem:cfs-signs}.
\end{remark}

\subsection{The all-split stratum and a Calabi--Yau threefold}\label{subsec:cfs}

This subsection makes Remark~\ref{rem:geometry}(iii) precise. We define a
threefold \(\mathcal R\) over \(\Z[1/2]\) whose \(\Fq\)-points are the classes of
all-split Richelot data with ordered points and ordered fixed points, and we show
that it is isomorphic to an open subscheme \(\mathcal Y^{\circ}\), defined by
explicit equations, of the model of Cynk, Freitag and Salvati
Manni \cite{CynkFreitagSalvatiManni2011}. By
\cite[\S3, Prop.~5]{CynkFreitagSalvatiManni2011} and the theta relations recalled
below, \(\mathcal Y^{\circ}(\mathbb C)\) is the part of the Siegel modular
threefold on which the ten fourth powers of the even theta constants do not
vanish; see
Remark~\ref{rem:cfs-model} for what this rests on. Section and proposition numbers of
\cite{CynkFreitagSalvatiManni2011} refer to its arXiv version. In this subsection
\(x_0,\dots,x_3\) and \(z\) denote the coordinates introduced below, not those
of Proposition~\ref{prop:hyp}.

To define \(\mathcal R\), put
\[
\Delta=\lambda(1-\lambda)\,ab\,(a-b)(1-a)(1-b)(\lambda-a)(\lambda-b),\qquad
\rho=(1-a)(1-b)(\lambda-a)(\lambda-b),
\]
and
\[
\mathcal R=\operatorname{Spec}\,\Z[\tfrac12][\lambda,a,b,\Delta^{-1}][s_1,s_2,s_3]
\big/\bigl(s_1^{2}-\lambda,\ s_2^{2}-ab,\ s_3^{2}-\rho\bigr).
\]
Since \(\lambda\), \(ab\), \(\rho\) and \(2\) are units, \(\mathcal R\) is finite
\'etale of degree \(8\) over the open subscheme \(\Delta\neq0\) of
\(\mathbb A^{3}_{\Z[1/2]}\). That subscheme is the moduli space of six ordered
distinct points on \(\PP^{1}\) in the normal form \((0,\infty,1,\lambda,a,b)\).
In particular \(\mathcal R\) is smooth over \(\Z[1/2]\).

For six distinct points \(P_1,\dots,P_6\) of \(\PP^{1}\) over a field, let
\(\iota_U\) be the involution of \(\PP^{1}\) that exchanges \(P_3,P_4\) and
exchanges \(P_5,P_6\). Let \(\iota_V\) exchange \(P_5,P_6\) and \(P_1,P_2\), and
let \(\iota_W\) exchange \(P_1,P_2\) and \(P_3,P_4\). Each of them is unique,
because a M\"obius transformation is determined by the images of three points.
Each exists, by the matrices in the proof of Lemma~\ref{lem:cfs-moduli}. If no
\(P_i\) is \(\infty\) and \(u,v,w\) are the monic quadratics with roots
\(\{P_1,P_2\}\), \(\{P_3,P_4\}\), \(\{P_5,P_6\}\), then \(\iota_U=\iota_{v,w}\).
Its fixed points are the roots of \(U=[v,w]\), together with \(\infty\) when
\(U\) has degree one. The same holds for \(V\) and \(W\).

\begin{lemma}\label{lem:cfs-moduli}
	Let \(k\) be a field of characteristic different from \(2\). Then \(\mathcal
	R(k)\) is in bijection with the set of \(\PGL_2(k)\)-orbits of pairs \((P,Q)\)
	of the following kind: \(P=(P_1,\dots,P_6)\) is a sequence of distinct
	points of \(\PP^{1}(k)\), and \(Q=(Q_1,\dots,Q_6)\) is a sequence of points of
	\(\PP^{1}(k)\) such that \(\{Q_1,Q_2\}\), \(\{Q_3,Q_4\}\), \(\{Q_5,Q_6\}\) are
	the sets of fixed points of \(\iota_U\), \(\iota_V\), \(\iota_W\). The group
	\(\PGL_2(k)\) acts freely on these pairs.
\end{lemma}

\begin{proof}
	The group \(\PGL_2(k)\) acts simply transitively on ordered triples of
	distinct points of \(\PP^{1}(k)\). So the action is free, and each orbit
	contains exactly one pair with \((P_1,P_2,P_3)=(0,\infty,1)\). For that pair
	put \((\lambda,a,b)=(P_4,P_5,P_6)\); then \(\Delta(\lambda,a,b)\neq0\). In
	this normal form the three involutions have the matrices
	\[
	A_W=\begin{pmatrix}0&\lambda\\1&0\end{pmatrix},\qquad
	A_V=\begin{pmatrix}0&ab\\1&0\end{pmatrix},\qquad
	A_U=\begin{pmatrix}-M_1&-M_0\\M_2&M_1\end{pmatrix},
	\]
	where \(M_0,M_1,M_2\) are the minors \eqref{eq:minors} of
	\(v=(x-1)(x-\lambda)\) and \(w=(x-a)(x-b)\). They satisfy
	\(A_W^{2}=\lambda I\), \(A_V^{2}=ab\,I\) and \(A_U^{2}=(M_1^{2}-M_0M_2)I=\rho I\),
	the last by Lemma~\ref{lem:minors}(iv).

	Let \(A\) be one of these matrices, with \(A^{2}=cI\). The fixed points of
	the involution are the eigenlines of \(A\), and its eigenvalues are the two
	square roots of \(c\). These are distinct, because \(c\neq0\) and the
	characteristic is not \(2\). A fixed point is \(k\)-rational exactly when its
	eigenvalue lies in \(k\). So an ordering of the fixed points by
	\(k\)-rational points is the same as a choice of \(s\in k\) with
	\(s^{2}=c\), the first fixed point being \(\ker(A-sI)\).

	Choose \(s_1,s_2,s_3\) in this way for \(\iota_W\), \(\iota_V\), \(\iota_U\).
	The point \((\lambda,a,b,s_1,s_2,s_3)\) of \(\mathcal R(k)\) determines the
	orbit, and every point of \(\mathcal R(k)\) arises in this way.
\end{proof}

When \(k=\Fq\), the pairs of the lemma are the data \(u,v,w\) with all six
quadratics \(u,v,w,U,V,W\) split, together with an ordering of all their roots.
Each triple of divisors of type \((e;r)=(1,1,1;1,1,1)\) in the set
\(\mathcal P\) of Proposition~\ref{prop:intrinsic} carries \(64\) such
orderings. Since \(\PGL_2(\Fq)\) acts freely on pairs, the orbit--stabiliser
relation gives \(\#\mathcal R(\Fq)=64\,M(1,1,1;1,1,1)\). In the normal form the
\(s_i\) are determined up to sign by \((\lambda,a,b)\), so
\begin{equation}\label{eq:Rcount}
	\#\mathcal R(\Fq)=\sum_{\lambda,a,b}
	\bigl(1+\psi_q(\lambda)\bigr)\bigl(1+\psi_q(ab)\bigr)\bigl(1+\psi_q(\rho)\bigr),
\end{equation}
the sum over \(\lambda,a,b\in\Fq\) with \(\Delta\neq0\). This is the sum of
Remark~\ref{rem:geometry}(ii).

The ordering of the six roots of \(f=uvw\) is a full level-two structure on
\(\Jac(C)\) for \(C\colon y^{2}=f\). The partition
\(\{P_1P_2|P_3P_4|P_5P_6\}\) gives the kernel of the Richelot isogeny. When
\(\dcoef\neq0\), the \(Q_i\) are the roots of \(UVW\), that is, the
Weierstrass points of the quotient curve. We use this reading only as
motivation.

We now recall the model of Cynk, Freitag and Salvati Manni. Let \(\Gamma=\Gamma_2[2]\cap\Gamma_{2,0}[4]\). We number the elements of
\((\Z/2)^{2}\) as \((b_0,b_1,b_2,b_3)=(00,01,10,11)\), as in
\cite[\S3, Prop.~5]{CynkFreitagSalvatiManni2011}. Put
\[
x_j=\vartheta\genfrac{[}{]}{0pt}{}{00}{b_j}^{2},\qquad
y_5=\prod_j\vartheta\genfrac{[}{]}{0pt}{}{00}{b_j},\qquad
y_4=-\vartheta\genfrac{[}{]}{0pt}{}{10}{01}^{4}-\vartheta\genfrac{[}{]}{0pt}{}{00}{11}^{4},
\]
of weights \(1,2,2\). By \cite[\S3, Prop.~5]{CynkFreitagSalvatiManni2011}, the
ring of modular forms of even weight for \(\Gamma\) is the even-weight part of
\(\mathbb C[y_5,x_0,x_1,x_2,x_3,y_4]\). The defining relations are
\cite[\S3, (1)]{CynkFreitagSalvatiManni2011}
\begin{equation}\label{eq:cfs1}
	y_5^{2}=x_0x_1x_2x_3,\qquad
	2y_5^{2}=x_0^{2}x_1^{2}+x_0^{2}x_3^{2}+x_1^{2}x_3^{2}
	+\bigl(x_0^{2}+x_1^{2}-x_2^{2}+x_3^{2}+y_4\bigr)y_4 .
\end{equation}
The variety \(\mathcal Y\subset\PP(1,1,1,1,2,2)\) defined by \eqref{eq:cfs1} is
the Satake compactification of \(\mathbb H_2/\Gamma\). The description of the
ring is derived in \cite{CynkFreitagSalvatiManni2011} by a sketched argument
(``going down process'').

In \cite[\S4]{CynkFreitagSalvatiManni2011} the coordinate \(y_4\) is replaced by
\[
z=2y_4+x_0^{2}+x_1^{2}+x_3^{2}-x_2^{2}
\]
(there again called \(y_4\)). Given the first relation, the second relation of
\eqref{eq:cfs1} becomes
\begin{equation}\label{eq:cfs2}
	y_5^{2}=x_0x_1x_2x_3,\qquad z^{2}=l_0l_1l_2l_3,
\end{equation}
with
\[
\begin{aligned}
	l_0&=x_0+x_1+x_2+x_3, & l_1&=x_0-x_1+x_2-x_3,\\
	l_2&=x_0+x_1-x_2-x_3, & l_3&=x_0-x_1-x_2+x_3,
\end{aligned}
\]
that is, \(l_j=\sum_i(-1)^{b_j\cdot b_i}x_i\). The relations \eqref{eq:cfs1}
and the formula for \(z\) depend on the numbering. \S2 of
\cite{CynkFreitagSalvatiManni2011} numbers the characteristics differently,
and in that numbering \eqref{eq:cfs1} does not hold. Every permutation of
\(x_0,\dots,x_3\) permutes the \(l_j\) up to an even number of signs, so
\eqref{eq:cfs2} and the open conditions below do not depend on the numbering.

The fourth powers of the ten even theta constants are expressed in these
coordinates as follows. Put \(X_j=\vartheta\genfrac{[}{]}{0pt}{}{b_j}{00}(2Z)\).
The relations \cite[p.~54]{vanGeemenNygaard1995} (see also
\cite[\S2]{CynkFreitagSalvatiManni2011}) read
\[
x_j=\sum_k(-1)^{b_j\cdot b_k}X_k^{2},\qquad\text{equivalently}\qquad l_j=4X_j^{2}.
\]
For the six even characteristics \(\genfrac{[}{]}{0pt}{}{c}{b}\) with
\(c\neq0\) we use the quadrics displayed before Theorem~1.8 of
\cite[p.~54]{vanGeemenNygaard1995}, which come from Riemann's theta formula:
\[
\vartheta\genfrac{[}{]}{0pt}{}{c}{b}(Z)^{2}=\sum_{d}(-1)^{b\cdot d}X_dX_{d+c},
\]
with \(X_d\) indexed by \(d\in(\Z/2)^{2}\). Since \(c\cdot b=0\), the terms for
\(d\) and \(d+c\) are equal. Hence
\[
\vartheta\genfrac{[}{]}{0pt}{}{c}{b}^{4}=4\bigl(X_dX_{d+c}\pm X_{d'}X_{d'+c}\bigr)^{2},
\]
where \(\{d,d+c,d',d'+c\}=(\Z/2)^{2}\) and the two even \(b\) give the two signs.

For example, \(\vartheta\genfrac{[}{]}{0pt}{}{10}{01}^{2}=2(X_0X_2-X_1X_3)\). Moreover
\(x_0^{2}+x_1^{2}-x_2^{2}-x_3^{2}=8(X_0^{2}X_2^{2}+X_1^{2}X_3^{2})\). With
\(\vartheta\genfrac{[}{]}{0pt}{}{00}{11}^{4}=x_3^{2}\), the definition of
\(y_4\) gives
\[
\begin{aligned}
z&=x_0^{2}+x_1^{2}-x_2^{2}-x_3^{2}-2\vartheta\genfrac{[}{]}{0pt}{}{10}{01}^{4}\\
&=8(X_0^{2}X_2^{2}+X_1^{2}X_3^{2})-8(X_0X_2-X_1X_3)^{2}=16X_0X_1X_2X_3 .
\end{aligned}
\]
This is the coordinate of the quotient map in
\cite[\S4]{CynkFreitagSalvatiManni2011}. In terms of \(x\) and \(z\) the six
fourth powers are therefore
\begin{equation}\label{eq:tpm}
	t^{\pm}_{\alpha\beta|\gamma\delta}=\tfrac14\bigl(l_\alpha l_\beta+l_\gamma l_\delta\bigr)\pm\tfrac12 z,
\end{equation}
where the partition \(\{\alpha\beta|\gamma\delta\}\) of \(\{0,1,2,3\}\) is
\(\{02|13\}\), \(\{01|23\}\), \(\{03|12\}\) for \(c=10,01,11\). For instance
\(\vartheta\genfrac{[}{]}{0pt}{}{10}{01}^{4}=t^{-}_{02|13}\). From the identities
\[
l_0l_1-l_2l_3=4(x_0x_2-x_1x_3),\quad
l_0l_2-l_1l_3=4(x_0x_1-x_2x_3),\quad
l_0l_3-l_1l_2=4(x_0x_3-x_1x_2)
\]
and \eqref{eq:cfs2} one gets
\(t^{+}_{\alpha\beta|\gamma\delta}\,t^{-}_{\alpha\beta|\gamma\delta}
=(l_\alpha l_\beta-l_\gamma l_\delta)^{2}/16\). For instance,
\begin{equation}\label{eq:tplus}
	2t^{+}_{03|12}=x_0^{2}-x_1^{2}-x_2^{2}+x_3^{2}+z,\qquad
	t^{+}_{03|12}\,t^{-}_{03|12}=(x_0x_3-x_1x_2)^{2}.
\end{equation}
Hence the ten fourth powers \(x_j^{2}\) and \(t^{\pm}\) are all nonzero exactly
on
\[
\mathcal Y^{\circ}\colon\quad
x_0x_1x_2x_3\,(x_0x_1-x_2x_3)(x_0x_2-x_1x_3)(x_0x_3-x_1x_2)\neq0 .
\]
We take this condition, together with \eqref{eq:cfs2}, as the definition of
\(\mathcal Y^{\circ}\) over \(\Z[1/2]\).

What rests on \cite{CynkFreitagSalvatiManni2011} and
\cite{vanGeemenNygaard1995} is only the identification of
\(\mathcal Y^{\circ}\) with the theta-open part of the Satake
compactification, namely Proposition~5 and the relations (1) of
\cite[\S3]{CynkFreitagSalvatiManni2011}, the relations between \(x_j\) and
\(X_k\), and the quadrics for \(\vartheta\genfrac{[}{]}{0pt}{}{c}{b}^{2}\).
The other identities above are elementary. Proposition~\ref{prop:cfs} and
Corollary~\ref{cor:cfs} concern \(\mathcal Y^{\circ}\) as defined by equations.

Since \(x_0\) is invertible on \(\mathcal Y^{\circ}\) and has weight one,
\(\mathcal Y^{\circ}\) lies in the chart \(x_0\neq0\) of
\(\PP(1,1,1,1,2,2)\). That chart is the affine space with coordinates
\(x_i/x_0\), \(z/x_0^{2}\), \(y_5/x_0^{2}\). So \(\mathcal Y^{\circ}\) is
affine, and \(\#\mathcal Y^{\circ}(\Fq)\) is the number of solutions of
\eqref{eq:cfs2} with \(x_0=1\) that satisfy the open condition.

On \(\mathcal R\) put \(m=a(1-a)(\lambda-b)\), a unit.

\begin{proposition}\label{prop:cfs}
	Consider the formulas
	\[
	\Phi\colon\ (x_0:x_1:x_2:x_3:z:y_5)=
	\bigl(m:\ s_2s_3:\ a\,s_1s_3:\ (1-a)(\lambda-b)\,s_1s_2:\ z_\Phi:\ m\,s_1s_2s_3\bigr)
	\]
	in \(\PP(1,1,1,1,2,2)\), where
	\(z_\Phi=a(a-1)(a-b)(b-\lambda)(\lambda-1)(ab-\lambda)\), and
	\[
	\Psi\colon\
	\begin{aligned}[t]
		\lambda&=\frac{x_2x_3}{x_0x_1},\qquad
		a=\frac{2x_2(x_0x_3-x_1x_2)}{x_1\bigl(x_0^{2}-x_1^{2}-x_2^{2}+x_3^{2}+z\bigr)},\qquad
		b=a\,\frac{x_1x_3}{x_0x_2},\\
		s_1&=\frac{y_5}{x_0x_1},\qquad
		s_2=\frac{a\,y_5}{x_0x_2},\qquad
		s_3=\frac{(1-a)(\lambda-b)\,y_5}{x_0x_3}.
	\end{aligned}
	\]
	They define mutually inverse isomorphisms
	\(\Phi\colon\mathcal R\to\mathcal Y^{\circ}\) and
	\(\Psi\colon\mathcal Y^{\circ}\to\mathcal R\) over \(\Z[1/2]\).
\end{proposition}

\begin{proof}
	Both schemes are affine. We show that \(\Phi\) and \(\Psi\) are well defined
	on coordinate rings and that both composites are the identity. The formulas
	for \(\Phi\) are weighted homogeneous with \(x_0=m\) a unit. So \(\Phi\) is
	the map into the chart \(x_0\neq0\) given by \(x_i/x_0\), \(z/x_0^{2}\),
	\(y_5/x_0^{2}\).

	(a) We show that \(\Phi\) maps \(\mathcal R\) to \(\mathcal Y^{\circ}\). By hand,
	\(x_0x_1x_2x_3=m^{2}s_1^{2}s_2^{2}s_3^{2}=y_5^{2}\) on \(\mathcal R\). The
	relation \(z_\Phi^{2}=l_0l_1l_2l_3\), with \(x_0,\dots,x_3\) as in \(\Phi\), is
	a polynomial identity in \(\Z[\tfrac12][\lambda,a,b,s_1,s_2,s_3]\) modulo
	\(s_1^{2}-\lambda\), \(s_2^{2}-ab\), \(s_3^{2}-\rho\). It was verified with
	computer algebra (item~(8) of Section~\ref{sec:comp}). The open conditions
	follow from
	\[
	\begin{aligned}
		x_0x_1x_2x_3&=m^{2}\,\lambda\,ab\,\rho,\\
		x_0x_2-x_1x_3&=a(a-1)(a-b)(b-\lambda)\,s_1s_3,\\
		x_0x_1-x_2x_3&=-a(a-1)(b-\lambda)(\lambda-1)\,s_2s_3,\\
		x_0x_3-x_1x_2&=-a(a-1)(a-b)(b-\lambda)(\lambda-1)\,s_1s_2 .
	\end{aligned}
	\]
	These are checked by hand using \(s_1^{2}=\lambda\), \(s_2^{2}=ab\),
	\(s_3^{2}=\rho\). For example,
	\(x_0x_2-x_1x_3=(1-a)(\lambda-b)s_1s_3(a^{2}-s_2^{2})\). The right-hand sides
	are units, since each \(s_i^{2}\) is a unit.

	(b) The map \(\Psi\) is regular on \(\mathcal Y^{\circ}\). Its denominators are
	\(x_0x_1\), \(x_0x_2\), \(x_0x_3\), and
	\(x_1(x_0^{2}-x_1^{2}-x_2^{2}+x_3^{2}+z)=2x_1t^{+}_{03|12}\). All are units;
	for the last one this follows from \eqref{eq:tplus}.

	(c) We show that \(\Psi\) maps \(\mathcal Y^{\circ}\) to \(\mathcal R\). By hand,
	using \(y_5^{2}=x_0x_1x_2x_3\), one gets \(s_1^{2}=x_2x_3/(x_0x_1)=\lambda\) and
	\(s_2^{2}=a^{2}x_1x_3/(x_0x_2)=ab\). The relation \(s_3^{2}=\rho\) means
	\((1-a)^{2}(\lambda-b)^{2}x_1x_2=\rho\,x_0x_3\) for the functions
	\(\lambda,a,b\) of \(\Psi\). It was verified with computer algebra modulo
	\eqref{eq:cfs2}.

	It remains to show that \(\Delta\) is a unit. The functions \(\lambda\),
	\(1-\lambda=(x_0x_1-x_2x_3)/(x_0x_1)\), \(a\),
	\(1-b/a=(x_0x_2-x_1x_3)/(x_0x_2)\) and \(b\) are visibly units, so
	\(a-b\) is a unit too. Put \(\varpi=x_0^{2}-x_1^{2}-x_2^{2}+x_3^{2}\),
	\(T=\varpi+z\) and \(A=2x_2(x_0x_3-x_1x_2)\), so that \(a=A/(x_1T)\). Then
	\[
	\begin{aligned}
		a-1&=\frac{A-x_1T}{x_1T}, &
		a-\lambda&=\frac{x_0A-x_2x_3T}{x_0x_1T},\\
		b-1&=\frac{x_3A-x_0x_2T}{x_0x_2T}, &
		b-\lambda&=\frac{x_1x_3A-x_2^{2}x_3T}{x_0x_1x_2T}.
	\end{aligned}
	\]
	For a numerator \(n=c_1-c_2T\) put \(n'=c_1-c_2(\varpi-z)\). Then
	\(nn'=(c_1-c_2\varpi)^{2}-c_2^{2}\,l_0l_1l_2l_3\), by \eqref{eq:cfs2}. Let
	\(\Pi=(x_0x_1-x_2x_3)(x_0x_2-x_1x_3)(x_0x_3-x_1x_2)\). The four products
	\(nn'\) are, respectively,
	\[
	-4\Pi,\qquad -4x_2^{2}\Pi,\qquad -4x_2^{2}\Pi,\qquad -4x_2^{2}x_3^{2}\Pi .
	\]
	These are polynomial identities in \(\Z[x_0,\dots,x_3]\), verified with
	computer algebra. The products are units on \(\mathcal Y^{\circ}\) and
	\(n'\) is regular, so each \(n\) is a unit. Hence \(a-1\), \(a-\lambda\),
	\(b-1\) and \(b-\lambda\) are units.

	(d) We show that \(\Psi\circ\Phi=\mathrm{id}\). By hand, \(\Phi\) pulls back
	\(x_2x_3/(x_0x_1)\) to \(\lambda\) and \(x_1x_3/(x_0x_2)\) to \(b/a\). Once
	\(a\) is known, it pulls back \(y_5/(x_0x_1)\), \(a\,y_5/(x_0x_2)\) and
	\((1-a)(\lambda-b)\,y_5/(x_0x_3)\) to \(s_1\), \(s_2\), \(s_3\). That \(\Phi\)
	pulls back the formula for \(a\) to \(a\) is an identity of rational functions
	on \(\mathcal R\), verified with computer algebra.

	(e) We show that \(\Phi\circ\Psi=\mathrm{id}\). By hand, using
	\(y_5^{2}=x_0x_1x_2x_3\), \(\Psi\) pulls back \(s_2s_3\), \(s_1s_3\),
	\(s_1s_2\) and \(s_1s_2s_3\) to \(m\,x_1/x_0\), \((1-a)(\lambda-b)\,x_2/x_0\),
	\(a\,x_3/x_0\) and \(m\,y_5/x_0^{2}\). So \(\Phi\circ\Psi\) sends
	\((x_0:x_1:x_2:x_3)\) to \((m/x_0)\cdot(x_0:x_1:x_2:x_3)\) and \(y_5\) to
	\((m/x_0)^{2}y_5\). For the coordinate \(z\), the identity
	\(z_\Phi=(m/x_0)^{2}\,z\) for the functions \(\lambda,a,b\) of \(\Psi\) was
	verified with computer algebra modulo \eqref{eq:cfs2}. Since \(m/x_0\) is a
	unit of weight zero, \(\Phi\circ\Psi\) is the identity of
	\(\mathcal Y^{\circ}\subset\PP(1,1,1,1,2,2)\).
\end{proof}

Five identities in this proof were verified only by machine: the relation
\(z_\Phi^{2}=l_0l_1l_2l_3\) in (a), the relation \(s_3^{2}=\rho\) after \(\Psi\)
in (c), the four products \(nn'\) in (c), the pull-back of \(a\) in (d), and the
identity for \(z\) in (e). Every other step is a short computation by hand, also
checked by machine; the scripts are listed in item~(8) of Section~\ref{sec:comp}.

\begin{corollary}\label{cor:cfs}
	Let \(\mathcal Y^{\circ}\) be the open subscheme of the model \eqref{eq:cfs2}
	of \cite{CynkFreitagSalvatiManni2011} over \(\Z[1/2]\) defined above. For
	every odd prime power \(q\),
	\[
	\#\mathcal Y^{\circ}(\Fq)=\#\mathcal R(\Fq)=q^{3}-12q^{2}+56q-91-t_q .
	\]
\end{corollary}

\begin{proof}
	The first equality holds by Proposition~\ref{prop:cfs}. We evaluate
	\eqref{eq:Rcount}; this also proves the two statements of
	Remark~\ref{rem:geometry}(ii) whose proof is deferred to this point.

	Fix \(\lambda\neq0,1\) and let \(\mathcal A=\Fq\setminus\{0,1,\lambda\}\).
	For a function \(h\) without zeros on \(\mathcal A\), the sum of
	\(\psi_q(h(a)h(b))\) over ordered pairs \(a\neq b\) in \(\mathcal A\) is
	\(\bigl(\sum_{\mathcal A}\psi_q(h)\bigr)^{2}-(q-3)\). We apply this to
	\(h=x\), \(h=(x-1)(x-\lambda)\) and \(h=x(x-1)(x-\lambda)\), for which
	\(\sum_{\mathcal A}\psi_q(h)\) equals \(-1-\psi_q(\lambda)\),
	\(-1-\psi_q(\lambda)\) and \(-a_\lambda(q)\), respectively. The second value
	holds because \(\sum_{x\in\Fq}\psi_q((x-1)(x-\lambda))=-1\).
	Since \(ab\rho=h(a)h(b)\) for the third \(h\), the sum over \(a\neq b\) in
	\(\mathcal A\) of \((1+\psi_q(ab))(1+\psi_q(\rho))\) is
	\[
	(q-3)(q-4)-3(q-3)+4\bigl(1+\psi_q(\lambda)\bigr)+a_\lambda(q)^{2}
	=(q-3)(q-7)+4\bigl(1+\psi_q(\lambda)\bigr)+a_\lambda(q)^{2}.
	\]

	Now multiply by \(1+\psi_q(\lambda)\) and sum over \(\lambda\neq0,1\). We use
	three facts: \(\sum_{\lambda\neq0,1}(1+\psi_q(\lambda))=q-3\);
	\((1+\psi_q(\lambda))^{2}=2(1+\psi_q(\lambda))\) for \(\lambda\neq0\); and the
	definition \eqref{eq:Tintro} of \(\cT(q)\). They give
	\[
	\#\mathcal R(\Fq)=(q-3)\bigl(q^{2}-10q+29\bigr)
	+\sum_{\lambda\neq0,1}a_\lambda(q)^{2}+\cT(q).
	\]

	It remains to evaluate the sum of the \(a_\lambda(q)^{2}\). Summing
	\(\psi_q((x-\lambda)(y-\lambda))\) over all \(\lambda\) gives \(q-1\) if
	\(x=y\) and \(-1\) otherwise. Hence
	\[
	\sum_{\lambda\in\Fq}a_\lambda(q)^{2}
	=q\sum_x\psi_q\bigl(x(x-1)\bigr)^{2}-\Bigl(\sum_x\psi_q\bigl(x(x-1)\bigr)\Bigr)^{2}
	=q(q-2)-1 .
	\]
	Also \(a_0(q)^{2}=a_1(q)^{2}=1\). So
	\(\sum_{\lambda\neq0,1}a_\lambda(q)^{2}=q^{2}-2q-3\). With
	\(\cT(q)=-t_q-q-1\) (Theorem~\ref{thm:T}(b)) the total is
	\(q^{3}-12q^{2}+56q-91-t_q\).
\end{proof}

The corollary also follows from \(\#\mathcal R(\Fq)=64\,M(1,1,1;1,1,1)\) and
Proposition~\ref{prop:intrinsic}(b).

\begin{remark}\label{rem:cfs-model}
	By \cite[\S3, Prop.~5]{CynkFreitagSalvatiManni2011} and the relations recalled
	above, \(\mathcal Y^{\circ}(\mathbb C)\) is the open part of the Siegel modular
	threefold \(\mathcal Y\) of \cite{CynkFreitagSalvatiManni2011}, the Satake
	compactification of \(\mathbb H_2/\Gamma\), on which the ten fourth powers of
	the even theta constants do not vanish. In this sense the cofactor in
	Corollary~\ref{cor:B}\textup{(a)} is the number of \(\Fq\)-points of the open
	part of that Siegel modular threefold. This reading rests on
	\cite[\S3, Prop.~5]{CynkFreitagSalvatiManni2011}, whose proof is sketched
	there, on the relations of \cite[p.~54]{vanGeemenNygaard1995} recalled above,
	and on the choice of \eqref{eq:cfs2} and the open condition as a model
	over \(\Z[1/2]\). Proposition~\ref{prop:cfs} and Corollary~\ref{cor:cfs} do
	not depend on it.
\end{remark}

\begin{remark}\label{rem:cfs-thomae}
	(i) Under \(\Phi\), each of the ten fourth powers \(x_j^{2}\), \(t^{\pm}\)
	becomes \(\pm m\) times one of the ten Thomae products of the sextuple
	\((0,\infty,1,\lambda,a,b)\) (factors involving \(\infty\) omitted); this
	was checked by machine and is not used elsewhere in the paper.

	(ii) Under \(\Phi\), \(z\) vanishes exactly where \(ab=\lambda\), that is,
	where \(\iota_V\) (\(x\mapsto ab/x\)) equals \(\iota_W\) (\(x\mapsto\lambda/x\)).
	There a single involution exchanges the points of each of the three pairs.
	For an affine representative in \(\Adm(\Fq)\) of the corresponding triple of
	divisors this is equivalent to \(\dcoef=0\). One
	direction is Proposition~\ref{prop:pencil}. Conversely, if one involution
	exchanges the roots of each of \(u,v,w\), these are fibres of the quotient
	map \(\PP^{1}\to\PP^{1}\) by the involution, so they lie in a pencil and
	\(\dcoef=0\). Thus the coordinate \(z=16X_0X_1X_2X_3\), a product of theta constants at
	\(2Z\), vanishes on \(\mathcal R\) exactly on the decomposable locus.
\end{remark}

The next lemma concerns the group \(\Gamma=\Gamma_2[2]\cap\Gamma_{2,0}[4]\). We
write \(M\in\operatorname{Sp}_4(\Z)\) as
\(M=\left(\begin{smallmatrix}A&B\\C&D\end{smallmatrix}\right)\) with
\(2\times2\) blocks.

\begin{lemma}\label{lem:cfs-group}
	For \(M\in\operatorname{Sp}_4(\Z)\) with \(C\equiv0\bmod2\) put
	\(M^{\sharp}=\left(\begin{smallmatrix}A&2B\\C/2&D\end{smallmatrix}\right)\).
	Then \(M^{\sharp}\in\operatorname{Sp}_4(\Z)\), \(M^{\sharp}(2Z)=2\,(M Z)\)
	for \(Z\in\mathbb H_2\), and
	\[
	\Gamma_2[2]\cap\Gamma_{2,0}[4]=\{M\in\Gamma_2[2]:M^{\sharp}\in\Gamma_2[2]\}.
	\]
\end{lemma}

\begin{proof}
	With \(\delta=\left(\begin{smallmatrix}2E&0\\0&E\end{smallmatrix}\right)\) one
	has \(M^{\sharp}=\delta M\delta^{-1}\). Since \(\delta\) multiplies the
	symplectic form by \(2\), \(M^{\sharp}\) is symplectic, and it is integral
	when \(C\equiv0\bmod 2\). Moreover
	\((A\cdot2Z+2B)(\tfrac12C\cdot2Z+D)^{-1}=2(AZ+B)(CZ+D)^{-1}\). For
	\(M\in\Gamma_2[2]\) one has \(A\equiv D\equiv E\) and \(2B\equiv0\bmod2\). So
	\(M^{\sharp}\in\Gamma_2[2]\) if and only if \(C/2\equiv0\bmod2\), that is,
	\(C\equiv0\bmod4\).
\end{proof}

Since \(Z\mapsto2Z\) comes from the isogeny \(v\mapsto2v\) from
\(\mathbb C^{2}/(\Z^{2}+Z\Z^{2})\) to \(\mathbb C^{2}/(\Z^{2}+2Z\Z^{2})\), the
group \(\Gamma\) can be read as full level two on both sides of this
\((2,2)\)-isogeny. The index of \(\Gamma\) in \(\Gamma_2[2]\) is \(8\), which
coincides with the degree of \(\mathcal R\) over the space of six ordered
points; we do not use or claim a moduli interpretation of \(\mathcal R\) in
terms of \(\Gamma\). The index is \(8\) because
\(M\mapsto C/2\bmod2\) is a surjective homomorphism from \(\Gamma_2[2]\) to the
symmetric \(2\times2\) matrices over \(\F_2\), with kernel \(\Gamma\). Its
values are symmetric because \(A^{t}C\) is symmetric and \(A\equiv E\); it is a
homomorphism because \(C_{MN}\equiv C_M+C_N\bmod4\); and it is surjective
because \(\left(\begin{smallmatrix}E&0\\2S&E\end{smallmatrix}\right)\in\Gamma_2[2]\)
for every symmetric \(S\).

\begin{remark}\label{rem:cfs-signs}
	We return to Remark~\ref{rem:geometry}(iv). By Theorem~\ref{thm:A},
	\(\sigma=(-1)^{j+1}r_1r_2r_3=-e_1e_2e_3\,r_1r_2r_3\). By
	Proposition~\ref{prop:intrinsic}(b),
	\(64\,M(e;r)=C_{e;r}(q)+\sigma t_q\) for every odd \(q\). Hence
	\[
	\sum_{e,r}e^{\varepsilon}r^{\rho}\cdot64\,M(e;r)
	=\sum_{e,r}e^{\varepsilon}r^{\rho}\,C_{e;r}(q)\;-\;64\,t_q\,
	\delta_{\varepsilon,111}\,\delta_{\rho,111}.
	\]
	So \(t_q\) occurs in exactly one character, the product
	\(\chi_0=e_1e_2e_3r_1r_2r_3\) of all six signs.

	The four characters \((\varepsilon;111)\) with \(|\varepsilon|\ge2\) of
	Theorem~\ref{thm:walsh} come from the factor \(q-5+2j=q-2-(e_1+e_2+e_3)\) in
	\eqref{eq:A}. This factor is the number of points of \(\PP^{1}\) outside
	the support of the triple. Its term \(-(e_1+e_2+e_3)\) moves \(\chi_0\) to
	the three characters with \(|\varepsilon|=2\). Indeed, inserting
	\eqref{eq:A} into the Walsh sums gives the \(t_q\)-part
	\(-q(q-1)(q-2)\,t_q\) for \(\varepsilon=111\) and \(q(q-1)\,t_q\) for
	\(|\varepsilon|=2\), as in Table~\ref{tab:walsh}.

	This suggests a statement about cohomology. The group \(H=(\Z/2)^{6}\) acts on \(\mathcal R\). Three generators change the
	sign of \(s_1\), \(s_2\) or \(s_3\); they reverse the order of \((Q_5,Q_6)\),
	\((Q_3,Q_4)\) or \((Q_1,Q_2)\). The other three are
	\[
	\begin{aligned}
		(\lambda,a,b,s_1,s_2,s_3)&\mapsto
		\bigl(\lambda^{-1},a^{-1},b^{-1},\,s_1/\lambda,\,s_2/(ab),\,s_3/(ab\lambda)\bigr),\\
		(\lambda,a,b,s_1,s_2,s_3)&\mapsto
		\bigl(\lambda^{-1},a/\lambda,b/\lambda,\,s_1/\lambda,\,s_2/\lambda,\,s_3/\lambda^{2}\bigr),\\
		(\lambda,a,b,s_1,s_2,s_3)&\mapsto(\lambda,b,a,s_1,s_2,s_3).
	\end{aligned}
	\]
	By Lemma~\ref{lem:cfs-moduli}, these exchange \(P_1,P_2\), resp. \(P_3,P_4\),
	resp. \(P_5,P_6\), followed by renormalisation with \(x\mapsto1/x\), resp.
	\(x\mapsto x/\lambda\), resp. the identity. Each of the six maps preserves
	the relations of \(\mathcal R\) and multiplies \(\Delta\) by a unit, and the
	six maps are commuting involutions.

	The action is not free. On the locus \(ab=\lambda\) the three involutions
	\(\iota_U,\iota_V,\iota_W\) coincide, and the product of the three
	exchanges of the \(P_i\) fixes these points. For \(h\in H\) let
	\(\mathcal R^{h}(\Fq)\) be the set of \(x\in\mathcal R(\overline{\F}_q)\) with
	\(\mathrm{Frob}_q(x)=h(x)\). We expect, from Lang's theorem for \(\PGL_2\),
	that \(\#\mathcal R^{h}(\Fq)=64\,M(e;r)\), where \(e\) and \(r\) record which
	generators occur in \(h\); we have not written this out. Granting it, the
	Grothendieck--Lefschetz trace formula, applied to the isotypic parts of
	\(H^{*}_{c}(\mathcal R_{\overline{\F}_q},\Ql)\), shows that \(t_q\) enters only
	the trace of Frobenius on the \(\chi_0\)-isotypic part.

	\emph{Question.} Let \(\ell\neq p\).
	\begin{enumerate}[label=\textup{(\alph*)},leftmargin=*]
		\item Is the semisimplification of
		\(H^{3}_{c}(\mathcal R_{\overline{\F}_p},\Ql)^{\chi_0}\) the sum of the
		two-dimensional representation attached to \(g\) and of Tate twists of
		\(\Ql\)?
		\item Are all other isotypic parts of
		\(H^{*}_{c}(\mathcal R_{\overline{\F}_p},\Ql)\) mixed Tate?
		\item Is the quotient of \(\mathcal R\) by \(\ker\chi_0\) birational to the
		double octic of Proposition~\ref{prop:hyp}?
	\end{enumerate}
\end{remark}

\section{The decomposable case and the 2-primary part}\label{sec:dec}

Let \((u,v,w)\) be a decomposable admissible triple whose class is a nonsquare.
By Proposition~\ref{prop:pencil}(ii), \(M_2(u,v)\neq0\), so
Proposition~\ref{prop:factors} applies with \(K=\F_{q^{2}}\): the fixed points
\(\alpha,\beta\in K\setminus\Fq\) of \(\iota\) are conjugate, \(E_1\colon Y^{2}=F(X)\)
is the curve \eqref{eq:E1}, and \(\#\Jac(C)(\Fq)=\#E_1(K)\). In this case
\(r_u=r_v=r_w=-1\), and the remaining datum is the triple of signs
\(e_u,e_v,e_w\).

\begin{proof}[Proof of Theorem~\ref{thm:C}]
	By \eqref{eq:E1},
	\[
	F(X)=c\prod_{s\in\{u,v,w\}}(X-x_s),\qquad c=u(\beta)v(\beta)w(\beta),\quad
	x_s=-\frac{s(\alpha)}{s(\beta)},
	\]
	so \(E_1\) has the three \(K\)-rational points \(P_s=(x_s,0)\) of order two. By
	the criterion of \(2\)-descent applied to the model
	\(Y'^{2}=\prod_s(X'-cx_s)\), the point \(P_s\) lies in \(2E_1(K)\) if and only if
	\(c(x_s-x_{s'})\) is a square in \(K\) for both \(s'\neq s\)
	\cite[Ch.~1, Thm.~4.1]{Husemoller2004}. Since \(-1\) is a square in \(K\), the
	condition for the pair \(\{u,v\}\) is symmetric, and
	\[
	c\,(x_u-x_v)=w(\beta)\bigl(u(\beta)v(\alpha)-u(\alpha)v(\beta)\bigr).
	\]
	Let \((M_2,M_1,M_0)\) be the minors of \((u,v)\) and \(R=\Res(u,v)\). By
	Lemma~\ref{lem:minors}(iii) and \eqref{eq:alphabeta},
	\begin{align*}
	u(\beta)v(\alpha)-u(\alpha)v(\beta)
	&=(\beta-\alpha)\bigl(M_2\alpha\beta+M_1(\alpha+\beta)+M_0\bigr)\\
	&=(\beta-\alpha)\,\frac{2(M_0M_2-M_1^{2})}{M_2}
	=-\frac{M_2}{2}(\beta-\alpha)^{3}.
	\end{align*}
	Every element of \(\Fq^{\times}\) is a square in \(K\), so
	\(c(x_u-x_v)\equiv w(\beta)(\beta-\alpha)\) modulo squares of \(K\). We compute
	the quadratic character of the two factors, using
	\(\psi_{q^{2}}(z)=\psi_q(z^{q+1})\).

	First, \(\vartheta=(\beta-\alpha)^{2}=4R/M_2^{2}\) is a nonsquare of \(\Fq\), so
	\[
	\psi_{q^{2}}(\beta-\alpha)=\vartheta^{(q^{2}-1)/4}
	=\bigl(\vartheta^{(q-1)/2}\bigr)^{(q+1)/2}=(-1)^{(q+1)/2}.
	\]
	Second, \(w(\beta)^{q+1}=w(\alpha)w(\beta)\), and
	\(M_2^{2}\,w(\alpha)w(\beta)=\prod_{w(\varrho)=0}W(\varrho)\). By
	Proposition~\ref{prop:pencil}(i), \(W\) is an \(\Fq^{\times}\)-multiple of
	\(V=w'u-wu'\), whose value at a root \(\varrho\) of \(w\) is
	\(w'(\varrho)u(\varrho)\); the scalar enters squared, and
	\(\prod_{\varrho}w'(\varrho)=-\Disc(w)\),
	\(\prod_\varrho u(\varrho)=\Res(w,u)\). Since \(\Res(w,u)\) is a nonsquare,
	\[
	\psi_{q^{2}}\bigl(w(\beta)\bigr)=\psi_q\bigl(-\Disc(w)\Res(w,u)\bigr)
	=-(-1)^{(q-1)/2}\,\psi_q(\Disc w).
	\]
	Multiplying,
	\(\psi_{q^{2}}\bigl(c(x_u-x_v)\bigr)=-(-1)^{q}\,\psi_q(\Disc w)=\psi_q(\Disc w)\).
	Thus \(c(x_u-x_v)\) is a square in \(K\) exactly when \(w\) splits over \(\Fq\),
	and cyclically. The criterion for \(P_u\) involves the pairs \(\{u,v\}\) and
	\(\{u,w\}\), hence the factors \(w\) and \(v\), which proves the first
	assertion.

	If \(k=3\), all three points of order two are divisible by two, which is
	\(E_1[4]\subseteq E_1(K)\). If \(k=2\), exactly one of them is, namely the one
	indexed by the irreducible factor; a group \(\Z/2^{a}\times\Z/2^{b}\) with
	\(1\le a\le b\) in which exactly one element of order two is divisible by two
	has \(a=1\) and \(b\ge2\). If \(k\le1\), none is, so there is no point of order
	four and the \(2\)-primary part is \(E_1[2]\).
\end{proof}

\begin{remark}\label{rem:two-part-data}
	The theorem concerns the group \(E_1(\F_{q^{2}})\), which has the same order as
	\(\Jac(C)(\Fq)\) but not the same structure. The group \(\Jac(C)(\Fq)[2]\) is
	generated by the classes of the Frobenius-stable pairs of roots and is
	\((\Z/2)^{4}\), \((\Z/2)^{3}\), \((\Z/2)^{2}\), \((\Z/2)^{2}\) for \(k=3,2,1,0\),
	so that \(16\mid N\) and \(8\mid N\) for \(k=3,2\) are already visible there. The
	content of the theorem is the structure of \(E_1(\F_{q^{2}})\) and the
	congruence \(N\equiv4\pmod8\) for \(k\le1\). The number \(k\) does not determine
	the \(2\)-adic valuation of \(N\) when \(k\ge2\): in the range of
	\S\ref{sec:comp} the exponent \(m\) of case (b) takes every value from \(2\) to
	\(6\). Case (a) occurs exactly when \(q\ge7\). Indeed, by Lemma~\ref{lem:Dq} the three
	quadratics lie on a line \(s_0=\mu s_1+\nu\) with \(\mu^{2}+\nu\) a nonsquare, the
	discriminant \(s_1^{2}-4\mu s_1-4\nu\) restricted to the line is an irreducible
	quadratic in \(s_1\), and the line carries \((q-1)/2\) split quadratics by
	Lemma~\ref{lem:B}(iii); any three of them form a triple with \(k=3\), and three
	exist when \(q\ge7\).
\end{remark}

\begin{remark}\label{rem:BL}
	For supersingular curves of this type the structure of the \(2\)-part is known
	\cite{CostelloKummer2018}. For ordinary ones, Bernstein and Lange construct
	Jacobians \(J\) over \(\F_p\) with rational \(2\)-torsion and
	\(\#J(\F_p)=16\ell\), \(\ell\) prime, isogenous to the Weil restriction of an
	elliptic curve \(E\) over \(\F_{p^{2}}\), and write that computer experiments
	suggest that the subgroup of order \(16\) of \(E(\F_{p^{2}})\) is
	\((\Z/4)\times(\Z/4)\) \cite[\S5]{BernsteinLange2014}. If \(C\) has an
	involution \(\iota\) of this kind and all six Weierstrass points of \(C\) are
	rational, then the orbits of \(\iota\) on the six roots form an admissible
	triple with \(\dcoef=0\) in which each factor splits; so \(k=3\), and Theorem~\ref{thm:C}(a) gives \(E_1[4]\subseteq E_1(K)\), hence
	the subgroup \((\Z/4)^{2}\) when \(16\) exactly divides the order. This applies
	to their observation when their curve \(E\) is the curve \(E_1\) of
	\eqref{eq:E1} or its conjugate; the shape of the \(2\)-part is not invariant
	under isogeny, and we have not compared the two normalisations beyond the
	construction.
\end{remark}

\section{Computations}\label{sec:comp}

The following computations accompany the proofs. They use exact integer and
polynomial arithmetic, and the scripts and logs are in Online Resource~1.
The proofs depend on them in two places. The first is the evaluation of the
\(20\) polynomials \(P_{e;r}\) in the proof of Theorem~\ref{thm:A}, a finite
computation with polynomials of degree at most four; it was done by computer
algebra (item~(3)), and, as explained in that proof, it also follows from the
enumeration of item~(4). The second is the set of polynomial identities named
after the proof of Proposition~\ref{prop:cfs} (item~(8)). The assembly of the
\(56\) entries of Table~\ref{tab:walsh} in the proof of
Proposition~\ref{prop:poly} consists of linear combinations of the displayed
values, which can be checked by hand; it was also carried out symbolically
(item~(3)).

In addition, the following statements made in remarks rest on machine checks
only: the range of the exponent in Remark~\ref{rem:two-part-data} (item~(7)),
the statement of Remark~\ref{rem:kummer} (item~(9)), Remark~\ref{rem:cfs-thomae}(i), and
the properties of the action of the group \(H\) stated in
Remark~\ref{rem:cfs-signs} (item~(8)). Nothing else in the paper uses these
statements.

\begin{enumerate}[label=\textup{(\arabic*)},leftmargin=*]
	\item \emph{Lemmas and intermediate sums.} The displayed identities in the
	proofs of Lemmas~\ref{lem:A}, \ref{lem:B}, \ref{lem:P1} and \ref{lem:P2} and of
	Proposition~\ref{prop:reduction} were tested on all instances over \(\Fq\) for
	\(q\in\{3,5,7,9,11,13\}\).
	\item \emph{Walsh sums.} All \(64\) sums \(\cS[\varepsilon;\rho]\) were compared
	with Table~\ref{tab:walsh}, with \(\cT(q)\) computed from \eqref{eq:Tintro},
	for every odd prime power \(q\le31\) by two independent programs.
	\item \emph{From Table~\ref{tab:walsh} to Table~\ref{tab:classes}.} The
	assembly of Proposition~\ref{prop:poly} from Lemmas~\ref{lem:P1}
	and~\ref{lem:P2}, and the inversion \eqref{eq:inversion} giving
	\(P_{e;r}=(q-5+2j)\,C_{e;r}\), were carried out symbolically in
	\(\Z[q,\cK]\).
	\item \emph{Class counts.} Each \(N(e;r)\) was compared with \eqref{eq:A} for
	\(q\in\{3,5,7,9,11,13\}\) by an enumeration of \(\Adm(\Fq)\) that decides the
	splitting of the six quadratics by counting their roots in \(\PP^{1}(\Fq)\),
	without characters or resultants.
	\item \emph{Projective form.} For \(q\in\{3,5,7,11\}\) the numbers of elements
	of \(\mathcal P\) of each type were enumerated; they equal
	\(q(q^{2}-1)\cdot\tfrac1{64}\bigl(C_{e;r}(q)+\sigma t_q\bigr)\) for all \(64\)
	types. This includes the cases \((q,j)=(5,0)\) and \((3,1)\), which are
	treated by hand in the proof of Proposition~\ref{prop:intrinsic}(b).
	\item \emph{The moment.} The identity \(\cK(q)=-t_q\) of Theorem~\ref{thm:T} was
	checked by two independent programs for all odd primes \(q\le599\) and for
	all \(36\) proper prime powers \(q\le6859\), among them \(3^{8}\), \(5^{5}\)
	and \(7^{4}\).
	\item \emph{Theorem~\ref{thm:C}.} The statement, including the group structure
	of the \(2\)-primary part, was compared with independently computed groups
	\(E_1(\F_{q^{2}})\) for every decomposable triple with nonsquare class over
	\(\Fq\), for all odd \(q\le17\).
	\item \emph{Subsection~\ref{subsec:cfs}.} The identities named after the proof
	of Proposition~\ref{prop:cfs} were verified symbolically, as were the
	elementary identities used in that proof and in the description of
	\(\mathcal Y^{\circ}\), the derivation of \(z=16X_0X_1X_2X_3\) and of the
	six \(t^{\pm}\), the statement of Remark~\ref{rem:cfs-thomae}(i), and the
	properties of the action of \(H\) stated in
	Remark~\ref{rem:cfs-signs}. As tests,
	\(\#\mathcal Y^{\circ}(\Fq)\), \(\#\mathcal R(\Fq)\) and
	\(q^{3}-12q^{2}+56q-91-t_q\) were compared for every odd prime power
	\(q\le31\), for \(q=37,41,49\) and for \(q=81,121,125,169\). For each of these \(q\le49\) other than \(q=5\) they were also compared
	with an enumeration of the admissible triples for which all six quadratics
	split.
	The map \(\Phi\) was checked to be injective, with
	\(\Psi\circ\Phi=\mathrm{id}\), on \(\mathcal R(\F_p)\) for \(p=7,13,31\).
	\item \emph{Remark~\ref{rem:kummer}.} The matrices \(W_P\) were computed for all
	admissible ordered triples over \(\Fq\) for \(q=3,5,7\), and for \(3000\)
	random ordered triples of separable monic quadratics for each of
	\(q=9,11,13\), of which the admissible ones were kept; the seed is recorded
	in Online Resource~1.
\end{enumerate}

\section*{Statements and Declarations}

\noindent\textbf{Competing interests.} The author has no competing interests
to declare that are relevant to the content of this article.

\smallskip
\noindent\textbf{Funding.} No funding was received for conducting this study.

\smallskip
\noindent\textbf{Data availability.} The data generated by the computations of Section~\ref{sec:comp} (logs and
outputs) are included in Online Resource~1, together with the scripts; no other
data were generated or analysed.

\smallskip
\noindent\textbf{Code availability.} The scripts and logs of the computations of
Section~\ref{sec:comp} are provided as Online Resource~1.

\smallskip
\noindent\textbf{Use of AI tools.} The research questions and the plan of the paper are the author's. The author used a large
language model (Claude, Anthropic) as an assistant: to write and run exploratory
and verification computations; to propose approaches, candidate proofs and
arguments, among which the author chose and which the author examined, corrected,
and accepted or rejected; to search the literature; and to draft and edit the
text. The author read every proof in the article, checked the statements,
computations and citations, and obtained reviews of the manuscript, some of them
generated with other AI systems, which the author directed and evaluated. The
author takes full responsibility for the content of the article. No AI tool is
an author of this work.

\bibliographystyle{amsplain}
\bibliography{references}

\end{document}